\documentclass[10pt]{amsart}
\usepackage[T1]{fontenc}
\usepackage{geometry}
\usepackage[latin1] {inputenc}
\usepackage{amsmath}
\usepackage{amsfonts, amssymb, textcomp}
\usepackage[colorlinks=flase, linkcolor=red,urlcolor=green, citecolor=blue]{hyperref}
\usepackage{subeqnarray}

\usepackage{latexsym}
\usepackage{fancyhdr}
\usepackage{longtable}
\usepackage{amsmath, amssymb}
\usepackage{graphicx}
\usepackage{enumerate}

\numberwithin{equation}{section}

\theoremstyle{plain}
\newtheorem{proposition}{Proposition}[section]
\newtheorem{theorem}[proposition]{Theorem}
\newtheorem{lemma}[proposition]{Lemma}
\newtheorem{corollary}[proposition]{Corollary}
\newtheorem{definition}[proposition]{Definition}
\newtheorem{example}[proposition]{Example}

\newtheorem{remark}[proposition]{Remark}

\newcommand{\RR}{\mathbb{R}}
\newcommand{\CC}{\mathbb{C}}
\newcommand{\NN}{\mathbb{N}}

\newcommand{\pr}{\operatorname{pr}}

\let\on=\operatorname

\newenvironment{demo}[1]{\par\smallskip\noindent{\bf #1.}}{\par\smallskip}

\makeatletter
\@namedef{subjclassname@2020}{%
  \textup{2020} Mathematics Subject Classification}
\makeatother

\newsavebox{\fmbox}
\newenvironment{fmpage}[1]
 {\begin{lrbox}{\fmbox}\begin{minipage}{#1}}
 {\end{minipage}\end{lrbox}\fbox{\usebox{\fmbox}}}

\title[On the regularization of sequences and associated weight functions]
{On the regularization of sequences and associated weight functions}

\author[G.~Schindl]{Gerhard Schindl}

\address{G.~Schindl: Fakult\"at f\"ur Mathematik, Universit\"at Wien, Oskar-Morgenstern-Platz~1, A-1090 Wien, Austria.}
\email{gerhard.schindl@univie.ac.at}

\begin{document}

\begin{abstract}
We revisit and generalize the geometric procedure of regularizing a sequence of real numbers with respect to a so-called regularizing function. This approach was studied by S. Mandelbrojt and becomes useful and necessary when working with corresponding classes of ultradifferentiable functions defined via weight sequences and analogous weighted spaces. In this note we also study non-standard situations for the construction yielding the (log-)convex minorant of a sequence and allow a ``blow-up'' for the regularizing function.
\end{abstract}

\thanks{This research was funded in whole by the Austrian Science Fund (FWF) project 10.55776/P33417}
\date{\today}

\maketitle

\section{Introduction}
Spaces of \emph{ultradifferentiable classes} are certain weighted subclasses of smooth functions; the weight is measuring the growth of the derivatives of the functions under consideration. Classically there do exist two, in general mutually distinct approaches: the more original method using a weight sequence $\mathbf{M}=(M_p)_{p\in\NN}$, see e.g. \cite{Komatsu73}, and the more recent technique using a weight function $\omega:[0,+\infty)\rightarrow[0,+\infty)$ in the sense of Braun-Meise-Taylor from \cite{BraunMeiseTaylor90}. Concerning these function spaces we refer also to \cite{BonetMeiseMelikhov07} and \cite{compositionpaper} and the citations therein. Moreover, analogous weighted structures are studied: We mention \emph{weighted spaces of sequences of complex numbers} or \emph{weighted spaces of formal power series,} see e.g. \cite{petzsche} and \cite{Borelmapalgebraity}, \emph{spaces of ultraholomorphic functions,} see e.g. \cite{Thilliezdivision} and \cite{optimalflat23}, generalized \emph{Gelfand-Shilov classes} see \cite{nuclearglobal2}, and \emph{weighted spaces of entire functions;} see \cite{weightedentireinclusion1} and \cite{weightedentireinclusion2} and in each instance also the references in these works. In all these settings growth and regularity assumptions on the weights are required and unavoidable and it turns out that crucial properties can be viewed as certain regularities of the weights. For a sequence $\mathbf{M}$ one of the most basic and important requirements is \emph{log-convexity:}
\begin{equation}\label{logconv}
\forall\;p\ge 1:\;\;\;M_p^2\le M_{p-1}M_{p+1},
\end{equation}
in the literature also denoted by $(M.1)$; see \cite{Komatsu73}. It turns out that this property is crucial for technical estimates in many proofs and, in particular, we mention that log-convexity implies that the ultradifferentiable classes are closed under the point-wise product; see e.g. \cite[Prop. 2.0.8]{diploma}.

A further useful fact is that log-convex sequences $\mathbf{M}$ can be expressed in terms of the corresponding \emph{associated function} $\omega_{\mathbf{M}}$ via \eqref{sequnderline} and \eqref{Komatsu73Prop32}; see Section \ref{assofunctsect} for more details about $\omega_{\mathbf{M}}$. This function appears frequently in both the weight sequence and weight function setting crucially. In the latter case $\omega_{\mathbf{M}}$ is also serving as a (counter-)example and is becoming relevant for comparison results between these settings; see \cite{BonetMeiseMelikhov07}.

The intimate relation between $\mathbf{M}$ and $\omega_{\mathbf{M}}$ expressed in terms of \eqref{sequnderline} resp. \eqref{Komatsu73Prop32} admits also applications in different contexts. We mention that it is unavoidable when proving the characterization of inclusion relations and stability properties in terms of the defining weights: For the ultradifferentiable setting see \cite[Prop. 4.6 \& Thm. 4.11]{compositionpaper}, for the weighted formal power series classes see \cite[Sect. 4.2]{Borelmapalgebraity}, for the Gelfand-Shilov spaces see \cite[Sect. 4 \& 6]{GelfandShilovincl}, and for the weighted entire setting we refer to \cite{weightedentireinclusion1} and \cite{weightedentireinclusion2}.\vspace{6pt}

Motivated by the importance of log-convexity one may ask:
\begin{itemize}
\item[$(*)$] When $\mathbf{M}$ violates \eqref{logconv}, is it possible to replace $\mathbf{M}$ by a log-convex sequence $\mathbf{L}$ such that the corresponding weighted function classes coincide? And how are the functions $\omega_{\mathbf{M}}$ and $\omega_{\mathbf{L}}$ related?

\item[$(*)$] More generally, when given an ``irregular behaving'' $\mathbf{M}$, can it be replaced by a ``more regular sequence'' without changing the corresponding weighted class?

\item[$(*)$] How can such a more regular sequence be obtained resp. computed from the given sequence $\mathbf{M}$?
\end{itemize}

S. Mandelbrojt treated these questions in \cite[Chapitre I]{mandelbrojtbook} by developing a purely geometric and abstract procedure: Consider a sequence $\mathbf{a}:=(a_p)_{p\in\NN}\in\RR^{\NN}$, formally one can even allow $a_p=+\infty$ for finitely many indices $p\ge 1$, and a so-called \emph{regularizing function} $\phi:[0,+\infty)\rightarrow[0,+\infty]$ which satisfies some basic growth properties; see Definition \ref{regufctdef}. Then construct the regularized sequence $\mathbf{a}^{\phi}=(a^{\phi}_p)_{p\in\NN}$, to be precise $(p,a^{\phi}_p)_{p\in\NN}$ is obtained from the given set $(p,a_p)_{p\in\NN}$. The crucial correspondence between $\mathbf{M}$ and $\mathbf{a}$ and their regularized counterparts is given by $a_p=\log(M_p)$, $M^{\phi}_p=\exp(a^{\phi}_p)$; i.e.
\begin{equation}\label{Mvsarelation}
\mathbf{a}=\log(\mathbf{M}),\hspace{30pt}\mathbf{M}^{\phi}=\exp(\mathbf{a}^{\phi}).
\end{equation}
This construction yields, as a ``by-product'', the notion of the so-called \emph{trace function.} The cases $\phi=+\infty$, i.e. $\phi(t)=+\infty$ for all $t\ge 0$ and which is formally not a regularizing function, and $\phi\neq +\infty$ were considered separately since the behavior of the regularized sequences is different. In the construction it turned out that for $\phi=+\infty$ the requirement $\lim_{p\rightarrow+\infty}\frac{a_p}{p}=+\infty$ avoids technical complications and via \eqref{Mvsarelation} this corresponds to $\lim_{p\rightarrow+\infty}(M_p)^{1/p}=+\infty$. Moreover, $\phi=+\infty$ yields the most prominent example the so-called \emph{log-convex minorant} $\mathbf{M}^{\on{lc}}$; see Section \ref{logminorantsect}. In the case $\phi\neq+\infty$ the example $\phi(t):=e^t$ is useful; this case is called \emph{``r\'{e}gularisation exponentielle''} in \cite[Chapitre I]{mandelbrojtbook}.\vspace{6pt}

Indeed, $\mathbf{M}^{\on{lc}}$ answers the first question above; see \cite[Thm. 2.15]{compositionpaper} and the citations there for ultradifferentiable classes and \cite[Thm. 5.7]{ultraholomstability} deals with the ultraholomorphic setting. It also turns out that $\omega_{\mathbf{M}}=\omega_{\mathbf{M}^{\on{lc}}}$ which is useful when working with this function. Concerning the second question we mention the classical work \cite{cartanmandelbrojt40} and \cite[Chapitre I, Sect. 1.9]{mandelbrojtbook}. In the decomposition result from \cite{Mandelbrojt40} the regularization w.r.t. $\phi(t)=e^t$ appeared crucially. Moreover, for $\mathbf{M}^{\on{lc}}$ the trace function is closely connected to $\omega_{\mathbf{M}}$, see \eqref{tildecomparison} and \eqref{Komatsu73Prop32}, and hence to a weight function in the sense of Braun-Meise-Taylor.\vspace{6pt}

The aim of this paper is to study this construction again in detail, to give a more complete picture in several directions and to extend Mandelbrojt's techniques:

\begin{itemize}
\item[$(A)$] When studying again \cite[Chapitre I]{mandelbrojtbook} in detail we found that the construction is formally not complete in the sense that the first term $a_0$ of $\mathbf{a}$ was not considered and that only non-negative slopes of the crucial straight lines in the construction were used; see the explanations on \cite[Chapitre I, p. 2]{mandelbrojtbook} and the normalization condition $\phi(0)\ge 1$ there. We close this technical gap and present a complete construction and in order to do so we are dealing now with a slightly more general notion for $\phi$ being a regularizing function; i.e. extending it to whole $\RR$, see Definition \ref{regufctdef}. Note that $a_0$ corresponds to $M_0$ and this value contains crucial information in the weighted settings.

\item[$(B)$] If $\phi=+\infty$, then study the situation when Mandelbrojt's basic assumption $\lim_{p\rightarrow+\infty}\frac{a_p}{p}=+\infty$ is violated.

\item[$(C)$] We focus in detail on the case when $\lim_{t\rightarrow T}\phi(T)=+\infty$ for some $T\in\RR$, i.e. allowing a ``blow-up-phenomenon'', since this was only mentioned briefly on \cite[Chapitre I, p. 2]{mandelbrojtbook}.
\end{itemize}

Indeed, these questions stem from very recent research. $(A)$ is motivated by \cite{anisolog} and \cite{GelfandShilovincl}: In these works the inclusion relations for Gelfand-Shilov classes are characterized. In order to proceed, for general \emph{anisotropic resp. multi-index sequences} $(M_{\alpha})_{\alpha\in\NN^d}$ in \cite{anisolog}, the definition and the construction of the log-convex minorant has to be generalized accordingly. It turned out that the higher-dimensional case requires a delicate induction argument and during the preparation of \cite{anisolog} we detected the gap mentioned in $(A)$.

$(B)$ was mainly inspired by \cite{ultraholomstability}: In the main results in \cite[Sect. 5]{ultraholomstability} log-convexity and the log-convex minorant of sequences are becoming crucial. In particular, in view of \cite[Thm. 5.8]{ultraholomstability} the behavior of $p\mapsto(M_p)^{1/p}$ and the failure of Mandelbrojt's basic assumption are interesting. Note that B. Rodr\'{i}guez-Salinas in \cite{RodriguezSalinas62} considered for ultraholomorphic classes also ``non-standard situations'' but focusing on the properties of the function classes and not giving an explicit construction for $\mathbf{M}^{\on{lc}}$; we refer to our detailed comments \cite[Rem. 5.5 \& 5.9]{ultraholomstability}.

Moreover, the following crucial fact is related to $(B)$: H. Komatsu used $\omega_{\mathbf{M}}$ in \cite{Komatsu73} but without assuming Mandelbrojt's basic assumption; i.e. not requiring necessarily $\lim_{p\rightarrow+\infty}(M_p)^{1/p}=+\infty$. But in this case $\omega_{\mathbf{M}}$ violates crucial growth requirements and, in particular, cannot be considered as weight function in the sense of Braun-Meise-Taylor.\vspace{6pt}

Summarizing, problems $(A)-(C)$ admit, on the one hand, a more detailed understanding and further applications in weighted settings and, on the other hand, can become relevant in general when studying the growth and regularity behavior of sequences of real numbers from an purely abstract point of view.\vspace{6pt}

The paper is structured as follows: First, in Section \ref{weightsequassofctsect} we introduce the notion of weight sequences $\mathbf{M}$ and study basic growth requirements for $\mathbf{M}$ and its consequences for $\omega_{\mathbf{M}}$. In particular, we study Komatsu's setting in Section \ref{Komatsusect} and recall and prove representation formulas for $\omega_{\mathbf{M}}$ under very general assumptions on $\mathbf{M}$. Section \ref{lcminorsect} is dedicated to the study of $(A)$ for $\phi=+\infty$ and hence we revisit the construction of the (log-)convex minorant $\mathbf{M}^{\on{lc}}$ of a given sequence $\mathbf{M}$ and provide a comparison with the results obtained by S. Mandelbrojt; see Section \ref{comparisonsection}. In Section \ref{nonstandardsect} we treat $(B)$ and study two ``degenerate/non-standard'' situations for the construction of $\mathbf{M}^{\on{lc}}$. The final Section \ref{generalsection} focuses on $(A)$ and $(C)$ for general regularizing functions and we give comparisons with the case $\phi=+\infty$. It turns out that $(C)$ and the non-standard case in Section \ref{case2sect} are similar.\vspace{6pt}

Finally, we wish to emphasize that this note is purely investigating the geometric regularization approach of $\mathbf{a}$ resp. $\mathbf{M}$ being related via \eqref{Mvsarelation} and the correspondence to the trace function and so to $\omega_{\mathbf{M}}$. It is known that more specific growth properties for $\mathbf{M}$ used in the theories of weighted settings transfer to properties for $\omega_{\mathbf{M}}$ and vice versa; see e.g. \cite[Lemma 2.4]{sectorialextensions}, \cite{subaddlike}, \cite{modgrowthstrange} and the references therein.\vspace{6pt}

\textbf{Acknowledgements.} The author wishes to thank the anonymous referee for the careful reading and efforts, the comments and the valuable suggestions which have improved and clarified the presentation of the results.

\section{Weight sequences and associated functions}\label{weightsequassofctsect}

\subsection{Basic notation}
We write $\NN:=\{0,1,2,\dots\}$ and $\NN_{>0}:=\{1,2,\dots\}$. Occasionally, we put $\RR_{>0}:=(0,+\infty)$ and $\RR_{\ge 0}:=[0,+\infty)$. Frequently we use the conventions $0^0:=1$, $\frac{1}{+\infty}=0$, $0\cdot(-\infty):=0$ and $p\cdot(-\infty):=-\infty$ for any $p\in\NN_{>0}$.

\subsection{Preliminaries}\label{preliminarysection}
Let $\mathbf{M}=(M_p)_{p\in\NN}\in\RR_{>0}^{\NN}$ be a given sequence of positive real numbers. Moreover, let us set $\mu_p:=\frac{M_p}{M_{p-1}}$ for $p\ge 1$ and $\mu_0:=1$. Hence we have $\frac{M_p}{M_0}=\mu_p\cdots\mu_1$ for all $p\in\NN$; the case $p=0$ gives the empty product. $\mathbf{M}$ is called \emph{normalized} if $1=M_0\le M_1$. For concrete applications and examples it is convenient to introduce the following class of sequences:
$$\hypertarget{LCset}{\mathcal{LC}}:=\{\mathbf{M}\in\RR_{>0}^{\NN}:\;\mathbf{M}\;\text{is normalized, log-convex},\;\lim_{p\rightarrow+\infty}(M_p)^{1/p}=+\infty\}.$$
\emph{Note:} For each $\mathbf{M}\in\hyperlink{LCset}{\mathcal{LC}}$ we have $\mu_0=1\le\mu_1\le\mu_2\le\dots$ and $\lim_{p\rightarrow+\infty}\mu_p=+\infty$. In particular $\mathbf{M}$ is non-decreasing and in fact there exists a one-to-one correspondence between sequences belonging to the class \hyperlink{LCset}{$\mathcal{LC}$} and sequences $\mu_0=1\le\mu_1\le\mu_2\le\dots$ and $\lim_{p\rightarrow+\infty}\mu_p=+\infty$ via setting $M_p:=\prod_{i=0}^p\mu_i$.

We write $\mathbf{M}\le\mathbf{N}$ if $M_p\le N_p$ for all $p\in\NN$ and $\mathbf{M},\mathbf{N}\in\RR_{>0}^{\NN}$ are called \emph{equivalent} if
$$\exists\;C\ge 1\;\forall\;p\in\NN:\;\;\;\frac{1}{C^{p+1}}N_p\le M_p\le C^{p+1}N_p;$$
i.e. if $0<\inf_{p\in\NN_{>0}}\left(\frac{M_p}{N_p}\right)^{1/p}\le\sup_{p\in\NN_{>0}}\left(\frac{M_p}{N_p}\right)^{1/p}<+\infty$. Equivalence preserves the corresponding weighted classes; therefore it is no restriction to assume always $M_0=1$ since otherwise $\mathbf{M}$ can be replaced by $\widetilde{\mathbf{M}}$ with $\widetilde{M}_p:=M_p/M_0$. Similarly, for each $\mathbf{M}\in\RR_{>0}^{\NN}$ we can find an equivalent sequence $\mathbf{N}$ which is even normalized.

\begin{itemize}
\item[$(*)$] Log-convexity for $\mathbf{M}$ is equivalent to the fact that $\mu=(\mu_p)_{p\ge 1}$ is non-decreasing, hence $\lim_{p\rightarrow+\infty}\mu_p$ exists in $\RR_{>0}\cup\{+\infty\}$.

\item[$(*)$] Moreover, it is known, see e.g. \cite[Lemma 2.0.4]{diploma}, that for log-convex $\mathbf{M}$ one has that $p\mapsto(M_p/M_0)^{1/p}$ is non-decreasing. Thus also $\lim_{p\rightarrow+\infty}(M_p)^{1/p}$ exists in $\RR_{>0}\cup\{+\infty\}$.

\item[$(*)$] Let $\mathbf{M}$ be log-convex, then
$$\forall\;p\in\NN_{>0}:\;\;\;\frac{M_p}{M_0}=\mu_1\cdots\mu_p\le\mu_p^p\Longrightarrow(M_p)^{1/p}\le(M_0)^{1/p}\mu_p,$$
and so $\lim_{p\rightarrow+\infty}(M_p)^{1/p}\le\lim_{p\rightarrow+\infty}\mu_p$. In particular, $\lim_{p\rightarrow+\infty}(M_p)^{1/p}=+\infty$ implies $\lim_{p\rightarrow+\infty}\mu_p=+\infty$ and the converse is valid for any $\mathbf{M}\in\RR_{>0}^{\NN}$; see e.g. \cite[p. 104]{compositionpaper}.
\end{itemize}

\begin{lemma}\label{preliminarysectlemma}
Let $\mathbf{M}\in\RR_{>0}^{\NN}$ and assume that $\lim_{p\rightarrow+\infty}\mu_p=:C<+\infty$ exists. Then we get $$\lim_{p\rightarrow+\infty}\mu_p=\lim_{p\rightarrow+\infty}(M_p)^{1/p}.$$
\end{lemma}
Summarizing, for any log-convex $\mathbf{M}$ one has $\lim_{p\rightarrow+\infty}\mu_p=\lim_{p\rightarrow+\infty}(M_p)^{1/p}$
and both limits exist in $\RR_{>0}\cup\{+\infty\}$.

\demo{Proof}
 By assumption
 $$\forall\;\epsilon>0\;\exists\;p_{\epsilon}\in\NN_{>0}\;\forall\;p> p_{\epsilon}:\;\;\;C-\epsilon\le\mu_p\le C+\epsilon\Leftrightarrow M_{p-1}(C-\epsilon)\le M_p\le(C+\epsilon)M_{p-1},$$
 hence iterating this gives $M_{p_{\epsilon}}(C-\epsilon)^i\le M_{p_{\epsilon}+i}\le(C+\epsilon)^iM_{p_{\epsilon}}$ and so $$\forall\;i\in\NN_{>0}:\;\;\;(M_{p_{\epsilon}})^{1/(p_{\epsilon}+i)}(C-\epsilon)^{i/(p_{\epsilon}+i)}\le(M_{p_{\epsilon}+i})^{1/(p_{\epsilon}+i)}\le(C+\epsilon)^{i/(p_{\epsilon}+i)}(M_{p_{\epsilon}})^{1/(p_{\epsilon}+i)}.$$ When $\epsilon>0$ is now arbitrary and fixed, then as $i\rightarrow+\infty$ we get $(M_{p_{\epsilon}})^{1/(p_{\epsilon}+i)}(C-\epsilon)^{i/(p_{\epsilon}+i)}\rightarrow C-\epsilon$ and $(C+\epsilon)^{i/(p_{\epsilon}+i)}(M_{p_{\epsilon}})^{1/(p_{\epsilon}+i)}\rightarrow C+\epsilon$. Finally, since $\epsilon>0$ can be taken arbitrary, $\lim_{p\rightarrow+\infty}(M_p)^{1/p}=C$ is verified.
\qed\enddemo

However, note that for Lemma \ref{preliminarysectlemma} the existence of the limits is crucial: Let $0<c<1<d<+\infty$ and set $M_0:=1$, $M_p:=c^p$ for all $p\ge 1$ odd and $M_p:=d^p$ for all $p\ge 1$ even and then
$$0=\liminf_{p\rightarrow+\infty}\mu_p<\liminf_{p\rightarrow+\infty}(M_p)^{1/p}=c<d=\limsup_{p\rightarrow+\infty}(M_p)^{1/p}<\limsup_{p\rightarrow+\infty}\mu_p=+\infty.$$

\subsection{The log-convex minorant}\label{logminorantsect}
Let $\mathbf{M}\in\RR_{>0}^{\NN}$ be given. The \emph{log-convex minorant} $\mathbf{M}^{\on{lc}}$ is defined to be the point-wise, i.e. w.r.t. relation $\le$, largest sequence among all log-convex sequences $\mathbf{L}$ satisfying $\mathbf{L}\le\mathbf{M}$.

So, if $\mathbf{L}$ is log-convex and $\mathbf{L}\le\mathbf{M}$, then already $\mathbf{L}\le\mathbf{M}^{\on{lc}}(\le\mathbf{M})$ is valid.

Obviously, $\mathbf{M}=\mathbf{M}^{\on{lc}}$ if and only if $\mathbf{M}$ is log-convex.

\subsection{Associated function}\label{assofunctsect}
Let $\mathbf{M}=(M_p)_{p\in\NN}\in\RR_{>0}^{\NN}$ be given. Then, the \emph{associated function} $\omega_{\mathbf{M}}: \RR_{\ge 0}\rightarrow\RR_{\ge 0}\cup\{+\infty\}$ is defined as follows:
\begin{equation}\label{assofunc}
\omega_{\mathbf{M}}(t):=\sup_{p\in \NN}\log\frac{M_0t^p}{M_p},\qquad t\ge 0,
\end{equation}
with the convention that $0^0:=1$. This ensures $\omega_{\mathbf{M}}(0)=0$ and $\omega_{\mathbf{M}}(t)\ge 0$ for any $t\ge 0$ since $\frac{t^0M_0}{M_0}=1$ for all $t\ge 0$. \eqref{assofunc} corresponds to \cite[$(3.1)$]{Komatsu73} and we immediately have that $\omega_{\mathbf{M}}$ is non-decreasing and satisfying $\lim_{t\rightarrow+\infty}\omega_{\mathbf{M}}(t)=+\infty$. In the literature $\omega_{\mathbf{M}}$ is often extended to $t\in\RR$ resp. even to $t\in\RR^d$ or $t\in\CC^d$ in a radial-symmetric way; i.e. replacing $t$ by $|t|$.\vspace{6pt}

Moreover, we set $\varphi_{\omega_{\mathbf{M}}}:=\omega_{\mathbf{M}}\circ\exp$ and define the \emph{Young-conjugate}
\begin{equation}\label{Youngconj}
\varphi^{*}_{\omega_{\mathbf{M}}}(t):=\sup_{s\in\RR}\{ts-\varphi_{\omega_{\mathbf{M}}}(s)\},\;\;\;t\ge 0.
\end{equation}
\emph{Note:} In the literature frequently this conjugate is defined by restricting to $s\ge 0$ and both definitions coincide if $\omega_{\mathbf{M}}$ is \emph{normalized;} i.e. if $\omega_{\mathbf{M}}(t)=0$ for all $t\in[0,1]$. In any case, $\varphi^{*}_{\omega_{\mathbf{M}}}$ is a convex function.

\subsection{On H. Komatsu's setting for $\omega_{\mathbf{M}}$}\label{Komatsusect}
In \cite[$(3.1)$]{Komatsu73} it was mentioned that
\begin{equation}\label{liminfcond}
\mathbf{M}_{\iota}:=\liminf_{p\rightarrow+\infty}(M_p/M_0)^{1/p}=\liminf_{p\rightarrow+\infty}(M_p)^{1/p}>0
\end{equation}
is crucial and basic for the definition of $\omega_{\mathbf{M}}$ according to \eqref{assofunc}: Under this assumption $\omega_{\mathbf{M}}$ is non-decreasing, continuous, increasing faster than $t\mapsto\log(t)^p$ for any $p\in\NN$, hence in particular $\lim_{t\rightarrow+\infty}\omega_{\mathbf{M}}(t)=+\infty$, and $\log(t)=o(\omega_{\mathbf{M}}(t))$ as $t\rightarrow+\infty$.

However, $\lim_{p\rightarrow+\infty}(M_p)^{1/p}=+\infty$ was not considered necessarily in \cite{Komatsu73} and, in general, \eqref{liminfcond} is \emph{not sufficient} to avoid extreme situations. The next result illustrates this fact and should be compared with Section \ref{case2sect}; in particular with Example \ref{nonstandardsectex} there.

\begin{lemma}\label{assofunctsectlemma}
Let $\mathbf{M}\in\RR_{>0}^{\NN}$ be given and set $\inf_{p\in\NN_{>0}}(M_p/M_0)^{1/p}=:\mathbf{M}_{\inf}\le\mathbf{M}_{\iota}$. Then we get:
\begin{itemize}
\item[$(i)$] $\mathbf{M}_{\iota}>0$ implies
$$\forall\;0\le t\le\mathbf{M}_{\inf}:\;\;\;\omega_{\mathbf{M}}(t)=0.$$
\item[$(ii)$] $\mathbf{M}_{\iota}<+\infty$ implies
$$\forall\;t>\mathbf{M}_{\iota}:\;\;\;\omega_{\mathbf{M}}(t)=+\infty.$$
\item[$(iii)$] If $\lim_{p\rightarrow+\infty}(M_p)^{1/p}=+\infty$, then $\omega_{\mathbf{M}}(t)<+\infty$ for all $t\ge 0$.
\end{itemize}
\end{lemma}
$(ii)$ and $(iii)$ together yield that $\lim_{p\rightarrow+\infty}(M_p)^{1/p}=+\infty$ is equivalent to $\omega_{\mathbf{M}}(t)<+\infty$ for any $t\ge 0$; see also \cite[Rem. 1]{nuclearglobal2} even in the \emph{anisotropic situation.}

\demo{Proof}
$(i)$ $\mathbf{M}_{\iota}>0$ holds if and only if $\mathbf{M}_{\inf}>0$ since $M_p>0$ for all $p\in\NN$. In this situation
$$\forall\;0<c\le\mathbf{M}_{\inf}\;\forall\;p\in\NN:\;\;\;M_p\ge M_0c^p,$$
which implies $\frac{M_0t^p}{M_p}\le\left(\frac{t}{c}\right)^p$ for all $p\in\NN$ and $t\ge 0$ and proves the desired property.

$(ii)$ Let $t>\mathbf{M}_{\iota}$ and choose $\epsilon>0$ such that $t>\epsilon>\mathbf{M}_{\iota}$. Then, by definition we can find infinitely many $p\in\NN_{>0}$ such that $\epsilon\ge(M_p/M_0)^{1/p}\Leftrightarrow\frac{M_0}{M_{p}}\ge\left(\frac{1}{\epsilon}\right)^{p}$ and so $\log\left(\frac{M_0t^p}{M_p}\right)\ge p\log\left(\frac{t}{\epsilon}\right)$ for infinitely many $p\in\NN_{>0}$. Hence this part is shown.

$(iii)$ $\lim_{p\rightarrow+\infty}(M_p)^{1/p}=+\infty$ gives
$$\forall\;t\ge 0\;\exists\;C\ge 1\;\forall\;p\in\NN:\;\;\;t^p\le CM_p\Leftrightarrow\log\left(\frac{t^p}{M_p}\right)\le\log(C),$$
hence by definition $\omega_{\mathbf{M}}(t)\le\log(M_0)+\log(C)=\log(M_0C)$ for all $t\ge 0$.
\qed\enddemo

\begin{remark}
\emph{In view of Lemma \ref{assofunctsectlemma} the ``extreme example'' $M_p:=d^p$, $d>0$, yields}
$$\omega_{\mathbf{M}}(t)=0,\;\;\;t\in[0,d],\hspace{30pt}\omega_{\mathbf{M}}(t)=+\infty,\;\;\;t\in(d,+\infty).$$

\emph{For the sake of completeness we mention that in \cite[Lemma 7.2]{solidassociatedweight} a ``converse extreme situation'' was considered: The function $t\mapsto\log(1+t)$ cannot be considered as an associated function for any $\mathbf{M}\in\RR_{>0}^{\NN}$; but when allowing $M_p=+\infty$ for infinitely many $p\in\NN$, then it is possible to treat this ``limiting function'' as well.}
\end{remark}

\subsection{Representation formulas for $\omega_{\mathbf{M}}$}
We recall now the explicit representation for $\omega_{\mathbf{M}}$ for more regular sequences; see also \cite[1.8 III, 1.8 V]{mandelbrojtbook} and \cite[Lemma 8.2.4, Prop. 8.2.5]{diploma} but there we have assumed $\mu_1\ge 1$.

\begin{lemma}\label{lemma1}
Let $\mathbf{M}=(M_p)_{p\in\NN}$ be log-convex and assume that $\lim_{p\rightarrow+\infty}(M_p)^{1/p}=+\infty$. Then we get
\begin{equation}\label{lemma1equ}
\omega_{\mathbf{M}}(t)=0,\;\text{for}\;t\in[0,\mu_1],\;\;\;\;\omega_{\mathbf{M}}(t)=\log\left(\frac{M_0t^p}{M_p}\right)\;\text{for}\;t\in[\mu_p,\mu_{p+1}],\;p\ge 1.
\end{equation}
In particular, $\omega_{\mathbf{M}}(\mu_p)=\log\left(\frac{M_0\mu_p^p}{M_p}\right)$ for all $p\ge 1$. This formula implies that $\omega_{\mathbf{M}}$ is continuous, non-decreasing and $\lim_{t\rightarrow+\infty}\omega_{\mathbf{M}}(t)=+\infty$. Moreover, $\omega_{\mathbf{M}}$ is increasing faster than $t\mapsto\log(t)^p$ for any $p\in\NN$ which gives $\log(t)=o(\omega_{\mathbf{M}}(t))$ as $t\rightarrow+\infty$. Finally, $t\mapsto\omega_{\mathbf{M}}(e^t)=\varphi_{\omega_{\mathbf{M}}}(t)$ is convex on $\RR$.
\end{lemma}

\demo{Proof}
First, by assumption recall that we get $\lim_{p\rightarrow+\infty}\mu_p=+\infty$ as well; see Lemma \ref{preliminarysectlemma} and Section \ref{preliminarysection}. Now let $t\ge 0$ be given such that $\mu_p\le t\le\mu_{p+1}$ for some $p\ge 1$. Then for any $q>p$ we have
$$\frac{M_0t^p}{M_p}\ge\frac{M_0t^q}{M_q}\Leftrightarrow\frac{M_q}{M_p}\ge t^{q-p}\Leftrightarrow\mu_{p+1}\cdots\mu_q\ge t^{q-p},$$
which is clearly satisfied since $p\mapsto\mu_p$ is non-decreasing by log-convexity. Similarly, if $q<p$, then
$$\frac{M_0t^p}{M_p}\ge\frac{M_0t^q}{M_q}\Leftrightarrow\frac{M_p}{M_q}\le t^{p-q}\Leftrightarrow\mu_{q+1}\cdots\mu_p\le t^{p-q},$$
which holds by analogous reasons. For $p=1$ we have $q=0$ and get $\frac{M_0t}{M_1}\ge 1\Leftrightarrow\frac{M_1}{M_0}=\mu_1\le t$. Note that for these arguments it is not necessary to assume that $\mu_{p+1}>\mu_p$.

Finally, we have $\omega_{\mathbf{M}}(\mu_1)=\log(\mu_1M_0/M_1)=\log(M_1M_0/(M_0M_1))=0$. Let $t\in[0,\mu_1]$, then $\frac{t^pM_0}{M_p}\le 1\Leftrightarrow t^p\le\frac{M_p}{M_0}=\mu_1\cdots\mu_p$ for all $p\ge 1$ and for $p=0$ we have $t^0/M_0=1$. Altogether this verifies that $\omega_{\mathbf{M}}$ vanishes on $[0,\mu_1]$. Alternatively, since by definition $\omega_{\mathbf{M}}$ is non-decreasing for this property it suffices to check $\omega_{\mathbf{M}}(\mu_1)=0$.
\qed\enddemo

By using this result we give a proof of the useful integral representation formula for $\omega_{\mathbf{M}}$ stated in \cite[$(3.11)$]{Komatsu73}; see \cite[1.8. III, 1.8 V]{mandelbrojtbook} and \cite[Prop. 8.2.5]{diploma} under slightly stronger normalization assumptions on the sequence.

\begin{lemma}\label{intformula}
Let $\mathbf{M}=(M_p)_{p\in\NN}\in\RR_{>0}^{\NN}$ be log-convex and such that $\lim_{p\rightarrow+\infty}(M_p)^{1/p}=+\infty$. We define the \emph{counting function} $\Sigma_{\mathbf{M}}$ by
\begin{equation*}\label{counting}
\Sigma_{\mathbf{M}}(t):=|\{p\ge 1:\;\;\;\mu_p\le t\}|,\;\;\;t\in[0,+\infty),
\end{equation*}
and get the following integral representation formula:
\begin{equation}\label{intrep}
\forall\;t\ge 0:\;\;\;\omega_{\mathbf{M}}(t)=\int_0^t\frac{\Sigma_{\mathbf{M}}(s)}{s}ds=\int_{\mu_1}^t\frac{\Sigma_{\mathbf{M}}(s)}{s}ds.
\end{equation}
\end{lemma}

\demo{Proof}
By definition of $\Sigma_{\mathbf{M}}$ it is clear that $\int_0^t\frac{\Sigma_{\mathbf{M}}(s)}{s}ds=0$ for all $0\le t<\mu_1$. For any $p\ge 1$ we get $$\int_{\mu_p}^{\mu_{p+1}}\frac{\Sigma_{\mathbf{M}}(s)}{s}ds=p\int_{\mu_p}^{\mu_{p+1}}\frac{1}{s}=p\log(\mu_{p+1}/\mu_p),$$
which equals $0$ if $\mu_p=\mu_{p+1}$. Thus, if $\mu_p\le t\le\mu_{p+1}$ with $p\ge 1$, then we have
\begin{align*}
\int_0^{t}\frac{\Sigma_{\mathbf{M}}(s)}{s}ds&=\int_{\mu_1}^{t}\frac{\Sigma_{\mathbf{M}}(s)}{s}ds=\sum_{q=1}^{p-1}q\log(\mu_{q+1}/\mu_q)+\int_{\mu_p}^t\frac{\Sigma_{\mathbf{M}}(s)}{s}ds=\sum_{q=1}^{p-1}q\log(\mu_{q+1}/\mu_q)+p\log(t/\mu_p)
\\&
=p\log(t)-\sum_{q=1}^p\log(\mu_q)=p\log(t)-\log(M_p/M_0)=\log\left(\frac{t^pM_0}{M_p}\right).
\end{align*}
By Lemma \ref{lemma1} we are done.
\qed\enddemo

Indeed, when inspecting the proofs of Lemmas \ref{lemma1} and \ref{intformula}, then by taking into account Lemma \ref{preliminarysectlemma} and the comments from Section \ref{preliminarysection} we get the following version:

\begin{lemma}\label{lemma1exotic}
Let $\mathbf{M}=(M_p)_{p\in\NN}$ be log-convex and assume that $\lim_{p\rightarrow+\infty}(M_p)^{1/p}=\lim_{p\rightarrow+\infty}\mu_p=:C<+\infty$. Then:
\begin{itemize}
\item[$(i)$] \eqref{lemma1equ} holds for all $0\le t<C$.

\item[$(ii)$] \eqref{intrep} is valid for all $0\le t<C$ and note that we have to restrict $\Sigma_{\mathbf{M}}$ to $[0,C)$.
\end{itemize}
\end{lemma}

\subsection{Associated weight sequence}\label{assosequsect}
Given $\mathbf{M}\in\RR_{>0}^{\NN}$ such that $\lim_{p\rightarrow+\infty}(M_p)^{1/p}=+\infty$ is valid let us introduce the sequence $\underline{\mathbf{M}}$ via
\begin{equation}\label{sequnderline}
\underline{M}_p:=M_0\sup_{t>0}\frac{t^p}{\exp(\omega_{\mathbf{M}}(t))},\;\;\;\forall\;p\in\NN.
\end{equation}
\emph{Note:} Since $\omega_{\mathbf{M}}(0)=0$, in \eqref{sequnderline} we can also write $\sup_{t\ge 0}$ and $\underline{M}_0=M_0$ is clear because $\omega_{\mathbf{M}}(t)\ge 0$, $\omega_{\mathbf{M}}(0)=0$, and so $\sup_{t>0}\frac{1}{\exp(\omega_{\mathbf{M}}(t))}=1$.\vspace{6pt}

We have the following properties for $\underline{\mathbf{M}}$:

\begin{itemize}
\item[$(i)$] $\underline{\mathbf{M}}$ is log-convex since
\begin{align*}
\sup_{t>0}\frac{t^p}{\exp(\omega_{\mathbf{M}}(t))}=\exp\left(\sup_{t>0}\{p\log(t)-\omega_{\mathbf{M}}(t)\}\right)=\exp\left(\sup_{s\in\RR}\{ps-\omega_{\mathbf{M}}(e^s)\}\right)=\exp(\varphi^{*}_{\mathbf{M}}(p));
\end{align*}
recall \eqref{Youngconj}.

\item[$(ii)$] We have
$$\forall\;p\in\NN:\;\;\;\underline{M}_p\le M_p,$$
since by definition:
$$\forall\;p\in\NN:\;\;\;M_0\sup_{t>0}\frac{t^p}{\exp(\omega_{\mathbf{M}}(t))}=M_0\sup_{t>0}\frac{t^p}{\sup_{q\in\NN}M_0t^q/M_q}\underbrace{\le}_{q=p}M_p.$$

\item[$(iii)$] If $\mathbf{M}$ is log-convex, then $\mathbf{M}\equiv\underline{\mathbf{M}}$, i.e.
$$\forall\;p\in\NN:\;\;\;\underline{M}_p\ge M_p,$$
because as seen above $\underline{M}_0=M_0$ is clear and by Lemma \ref{lemma1} we get for any $p\in\NN_{>0}$ that
$$M_0\sup_{t>0}\frac{t^p}{\exp(\omega_{\mathbf{M}}(t))}\ge M_0\frac{\mu_p^p}{\exp(\omega_{\mathbf{M}}(\mu_p))}=M_p.$$
\end{itemize}

Gathering, by $(i)$ and $(ii)$ we immediately see that $\underline{\mathbf{M}}\le\mathbf{M}^{\on{lc}}$ and in the next section we prove that in fact equality holds.

Analogous statements hold if $\lim_{p\rightarrow+\infty}(M_p)^{1/p}=C\in(0,+\infty)$ instead: Then in \eqref{sequnderline} we have to replace $\sup_{t>0}$ by $\sup_{0<t<C}$ and use Lemma \ref{lemma1exotic} instead of Lemma \ref{lemma1} in $(iii)$.

\section{Constructing the (log-)convex minorant}\label{lcminorsect}
A sequence $\mathbf{a}:=(a_p)_{p\in\NN}\in\RR^{\NN}$ is called \emph{convex,} if the epigraph of the set $\{(p,a_p): p\in\NN\}$ is a convex set in $\RR^2$.

The \emph{convex minorant} $\mathbf{a}^c$ of $\mathbf{a}$ is defined to be the largest sequence w.r.t. relation $\le$, i.e. point-wise, among all convex sequences $\mathbf{b}\le\mathbf{a}$.

\subsection{Geometric construction of the convex minorant}\label{normgeomconstrsect}
Let $\mathbf{a}:=(a_p)_{p\in\NN}\in\RR^{\NN}$ be given and assume in this section that
\begin{equation}\label{basicregprop}
\lim_{p\rightarrow+\infty}\frac{a_p}{p}=+\infty.
\end{equation}
Moreover let us set $S_p:=(p,a_p)$ for $p\in\NN$ and $\mathcal{S}:=\{S_p: p\in\NN\}$.

Note:

\begin{itemize}
\item[$(a)$] The forthcoming geometric construction is not changing when we formally allow $a_p=+\infty$ for \emph{finitely many} $p\ge 1$.

\item[$(b)$] On the other hand, $\lim_{p\rightarrow+\infty}\frac{a_p}{p}=+\infty$ implies $a_p\ge p$ for all $p$ sufficiently large and hence, in particular, that
    \begin{equation}\label{basicregprop1}
    \liminf_{p\rightarrow+\infty}a_p>-\infty.
    \end{equation}
    In any case, \eqref{basicregprop1} implies
    \begin{equation}\label{basicregprop2}
\liminf_{p\rightarrow+\infty}\frac{a_p}{p}>-\infty,
\end{equation}
because $a_p\le\frac{a_p}{p}\Leftrightarrow 1\le p$ for all $p\ge 1$ such that $a_p\le 0$ and $\frac{a_p}{p}\ge 0$ for all $p\ge 1$ such that $a_p>0$.

\item[$(c)$] The value $\frac{a_p}{p}$ is precisely the slope of the straight line connecting the origin $(0,0)$ and the point $S_p$, $p\ge 1$.
\end{itemize}

In order to construct the convex minorant of $\mathbf{a}$ the idea is to project down vertically each $a_p$, $p\ge 1$, to a set of straight lines $\{D_k: k\in\RR\}$ called ``supporting lines''. Here $D_k$ is a straight line having slope $k\in\RR$ and is given by the formula
$$D_k: t\mapsto tk+d_k,\hspace{15pt}k\in\RR.$$

Formally, we can also consider the slope $k=-\infty$ with $D_{-\infty}$ denoting the $y$-axis $\{(0,y): y\in\RR\}$. $D_k$ are obtained iteratively as follows:\vspace{6pt}

\begin{itemize}
\item[$(1)$] Let the starting point $S_0=(0,a_0)$ be given and consider all straight lines passing through this point with slope $k$. For sufficiently small values $k>-\infty$ all points $S_p$, $p\ge 1$, will lie strictly above these lines because we have \eqref{basicregprop2}; see comment $(b)$ before. So the only point belonging to $\mathcal{S}$ and lying on these lines is $S_0$ and the ``limiting case'' $k=-\infty$ corresponds to the $y$-axis $\{(0,y): y\in\RR\}$. Then we let $k$ increase and we find some minimal slope $k_1$ and a minimal index $p_1\ge 1$ such that $S_{p_1}$ lies on the straight line passing through $S_0$ and with slope $k_1$.

    We denote this line by $D_{k_1}$ and its slope is given by $k_1=\frac{a_{p_1}-a_0}{p_1}$ which might be negative but $k_1>-\infty$. Thus we get
    $$D_{k_1}: t\mapsto\frac{a_{p_1}-a_0}{p_1}t+a_0.$$
    Let us set $p_0:=0$ and let $D_k$, with $-\infty<k\le k_1$, be the lines passing through $S_0$ with slope $k$.

\item[$(2)$] Then consider all straight lines passing through $S_{p_1}$ and having slope $k\ge k_1$. When increasing $k$ we can find, similarly as before, a minimal slope $k_2\ge k_1$ and a minimal index $p_2\ge p_1+1$ such that $S_{p_2}$ lies on the straight line passing through $S_{p_1}$ and having slope $k_2$.

    This line is denoted by $D_{k_2}$, its slope is given by $k_2=\frac{a_{p_2}-a_{p_1}}{p_2-p_1}$, which might be negative again, but it is well-defined since $p_2\ge p_1+1$. Note that $k_2=k_1$ might be valid and in this case $S_{p_2}$ already lies on $D_{k_1}$. Similarly, for $k_1<k<k_2$ let $D_k$ be the lines passing through $S_{p_1}$ and having slope $k$. When $k_1=k_2$ this set of lines is empty and $D_{k_1}=D_{k_2}$.

\item[$(3)$] Then proceed by iteration, so when given a point $S_{p_i}$ we can find a point $S_{p_{i+1}}$ with $p_{i+1}\ge p_i+1$ and a straight line $D_{k_{i+1}}$ with slope $k_{i+1}\ge k_i$ and passing through the points $S_{p_i}$ and $S_{p_{i+1}}$. $D_k$ with $k_i<k<k_{i+1}$ are defined to be the lines with slope $k$ and passing through $S_{p_i}$.

\item[$(4)$] So $D_k: t\mapsto kt+d_k$ and the lines $D_{k_i}$ are given by $t\mapsto tk_i+d_{k_i}$ with
$$k_i=\frac{a_{p_i}-a_{p_{i-1}}}{p_{i}-p_{i-1}},\hspace{30pt}d_{k_i}=\frac{p_ia_{p_{i-1}}-p_{i-1}a_{p_i}}{p_i-p_{i-1}},\;\;\;i\ge 1,$$
because $$d_{k_i}=a_{p_i}-\frac{a_{p_i}-a_{p_{i-1}}}{p_{i}-p_{i-1}}p_i=\frac{p_ia_{p_i}-p_{i-1}a_{p_i}-p_ia_{p_i}+p_ia_{p_{i-1}}}{p_i-p_{i-1}}=\frac{p_ia_{p_{i-1}}-p_{i-1}a_{p_i}}{p_i-p_{i-1}}.$$
As mentioned before one might have $k_{i+1}=k_i$ but $\lim_{i\rightarrow+\infty}k_{i}=\lim_{i\rightarrow+\infty}p_{i}=+\infty$ since $\lim_{p\rightarrow+\infty}\frac{a_p}{p}=+\infty$; see again assumption \eqref{basicregprop}. Moreover, for sufficiently large $i$ we have that $d_{k_i}<0$.
\end{itemize}

Then consider the set $\widetilde{S}=(\widetilde{S}_p)_{p\in\NN}$ with $\widetilde{S}_p=(p,\widetilde{a}_p)$ and the $y$-coordinate of each point is obtained by projecting the points $S_p$ vertically down onto the set of straight lines $D_k$. This yields a sequence $\widetilde{\mathbf{a}}=(\widetilde{a}_p)_{p\in\NN}$ and we get:
\begin{itemize}
\item[$(*)$] $\widetilde{a}_{p_i}=a_{p_i}$ for all $i\in\NN$,

\item[$(*)$] and for $p\ge 1$ with $p_i<p<p_{i+1}$:
$$\widetilde{a}_p=pk_i+d_{k_i}=\frac{a_{p_i}-a_{p_{i-1}}}{p_{i}-p_{i-1}}p+\frac{p_ia_{p_{i-1}}-p_{i-1}a_{p_i}}{p_i-p_{i-1}}.$$
\end{itemize}
In other words, by the geometric construction we have
\begin{equation}\label{regsequformula}
\forall\;p\in\NN:\;\;\;\widetilde{a}_p=\sup_{k\in\RR}\{kp+d_k\}=\sup_{k\in\RR}\{kp-A(k)\},
\end{equation}
with
\begin{equation}\label{tracefunction}
A(k):=-d_k.
\end{equation}
$A:\RR\rightarrow\RR$ is called the \emph{trace function} and $-A(k)$ is precisely the $y$-coordinate of the unique intersection point of $D_k$ with $\{(0,y): y\in\RR\}$; i.e. with the $y$-axis. Obviously $A$ is continuous and non-decreasing and we put as limiting case $A(-\infty):=-a_0$ since $k=-\infty$ corresponds to the $y$-axis itself and $S_0$ is the only point lying on this line. With the conventions $0\cdot(-\infty):=0$ and $p\cdot(-\infty):=-\infty$ for any $p\ge 1$ we get that in \eqref{regsequformula} we can consider $k\in\{-\infty\}\cup\RR$.

Note that from \eqref{regsequformula} it follows that
\begin{itemize}
\item[$(*)$] $\widetilde{\mathbf{a}}\equiv\mathbf{a}^c$,

\item[$(*)$] $\lim_{p\rightarrow+\infty}\frac{\widetilde{a}_p}{p}=+\infty$ and

\item[$(*)$] $\widetilde{a}_0=\sup_{k\in\RR}\{d_k\}=a_0$ since each $D_k$ with $-\infty<k\le k_1$ intersects the $y$-axis at $S_0$.

\item[$(*)$] We have
$$\forall\;i\in\NN:\;\;\;\widetilde{a}_{p_i}=a^c_{p_i}=a_{p_i},$$
and so the regularized sequence coincides with the original given one for infinitely many indices since $\lim_{i\rightarrow+\infty}p_i=+\infty$ by basic assumption \eqref{basicregprop}. In particular, we always have $\widetilde{a}_0=a^c_0=a_0$.
\end{itemize}

For any fixed $k\in\RR$ we see that $d_k$ can be expressed in terms of the given sequence $\mathbf{a}$ by
\begin{equation}\label{dkequ}
d_k=\inf_{p\in\NN}\{a_p-pk\}.
\end{equation}
This holds since the line passing through $S_p$ and having slope $k$ intersects the $y$-axis at the point $(0,a_p-pk)$ and by construction and definition of the lines $D_k$ we have \eqref{dkequ}. For $p=0$ we get as intersection point $(0,a_0)$ and clearly all lines passing through $S_0$ intersect the $y$-axis at $S_0$. Thus one has
\begin{equation}\label{trace}
\forall\;k\in\RR:\;\;\;A(k)=-d_k=-\inf_{p\in\NN}\{a_p-pk\}=\sup_{p\in\NN}\{pk-a_p\}.
\end{equation}
\emph{Note:} The conventions $0\cdot(-\infty):=0$ and $p\cdot(-\infty):=-\infty$ for any $p\ge 1$ give that \eqref{trace} even holds for $k=-\infty$ since then $A(-\infty)=-a_0$.

\subsection{Obtaining the log-convex minorant}\label{lcminorantconstrsect}
Let now $\mathbf{M}=(M_p)_{p\in\NN}\in\RR_{>0}^{\NN}$ be given and note that in applications for weighted structures the sequence $\mathbf{M}$ is the crucial one. The idea is to apply the geometric construction from the previous Section \ref{normgeomconstrsect} to
\begin{equation}\label{logconvtrafo}
a_p:=\log(M_p),
\end{equation}
and so in order to guarantee \eqref{basicregprop} we have to assume $\lim_{p\rightarrow+\infty}(M_p)^{1/p}=+\infty$. Note that formally we can also allow $M_p=+\infty$ for finitely many indices $p\ge 1$.

The sequence $\mathbf{M}^{\on{lc}}$ is then given by
\begin{equation}\label{logconvtrafo1}
M^{\on{lc}}_p=\exp(\widetilde{a}_p),
\end{equation}
so the log-convex minorant is obtained by considering $\mathbf{M}\mapsto\mathbf{a}\mapsto\widetilde{\mathbf{a}}\mapsto\mathbf{M}^{\on{lc}}$ via \eqref{logconvtrafo}, \eqref{logconvtrafo1} and the regularization from the previous Section. First, by \eqref{trace} we get for any $k\in\RR$:
$$A(k)=\sup_{p\in\NN}\{pk-a_p\}=\sup_{p\in\NN}\{pk-\log(M_p)\}=\sup_{p\in\NN}\log\left(\frac{e^{pk}}{M_p}\right)=:\widetilde{\omega}_{\mathbf{M}}(e^k)(=\varphi_{\widetilde{\omega}_{\mathbf{M}}}(k)),$$
i.e.
\begin{equation}\label{Avsomega}
A=\widetilde{\omega}_{\mathbf{M}}\circ\exp,
\end{equation}
and \eqref{logconvtrafo1} and \eqref{regsequformula} yield
\begin{align*}
\forall\;p\in\NN:\;\;\;M^{\on{lc}}_p&=\exp(\widetilde{a}_p)=\exp\left(\sup_{k\in\RR}\{kp-A(k)\}\right)=\sup_{k\in\RR}\frac{e^{kp}}{\exp(A(k))}=\sup_{k\in\RR}\frac{e^{kp}}{\exp(\widetilde{\omega}_{\mathbf{M}}(e^k))}
\\&
=\sup_{s>0}\frac{s^p}{\exp(\widetilde{\omega}_{\mathbf{M}}(s))}.
\end{align*}
For $p=0$ we point out that $M^{\on{lc}}_0=M_0=\sup_{s>0}\frac{1}{\exp(\widetilde{\omega}_{\mathbf{M}}(s))}$ since $\widetilde{\omega}_{\mathbf{M}}$ is by definition non-decreasing and continuous and $\lim_{s\rightarrow 0}\widetilde{\omega}_{\mathbf{M}}(s)=\lim_{t\rightarrow-\infty}A(t)=-a_0=-\log(M_0)$. Moreover, we can consider in the supremum all values $s\ge 0$.\vspace{6pt}

Let us now compare the sequence $\underline{\mathbf{M}}$ defined in \eqref{sequnderline} with $\mathbf{M}^{\on{lc}}$:

\begin{itemize}
\item[$(*)$] If $M_0=1$, i.e. $a_0=0$, then $\underline{\mathbf{M}}$ coincides with $\mathbf{M}^{\on{lc}}$ because $\omega_{\mathbf{M}}\equiv\widetilde{\omega}_{\mathbf{M}}$.

\item[$(*)$] If $M_0\neq 1$, then for sufficiently small values $s\ge 0$ one might have $\widetilde{\omega}_{\mathbf{M}}(s)<0$. More precisely, by definition $\widetilde{\omega}_{\mathbf{M}}(s)\ge 0$ for all $s\ge 0$ if and only if $a_0=\log(M_0)\le 0$.

    Thus, if $M_0>1$ then we see that $\widetilde{\omega}_{\mathbf{M}}(s)$ is negative for all sufficiently small $s\ge 0$ and hence, in particular, $\widetilde{\omega}_{\mathbf{M}}$ is not a weight function in the sense of Braun-Meise-Taylor; see \cite{BraunMeiseTaylor90}.

\item[$(*)$] However, by recalling \eqref{assofunc} in any case we have
    \begin{equation}\label{tildecomparison}
    \forall\;s\ge 0:\;\;\;\omega_{\mathbf{M}}(s)=\log(M_0)+\widetilde{\omega}_{\mathbf{M}}(s),
    \end{equation}
    and so
    \begin{equation}\label{Komatsu73Prop32}
    \forall\;p\in\NN:\;\;\;M^{\on{lc}}_p=\sup_{s>0}\frac{s^p}{\exp(\widetilde{\omega}_{\mathbf{M}}(s))}=M_0\sup_{s>0}\frac{s^p}{\exp(\omega_{\mathbf{M}}(s))}=\underline{M}_p;
\end{equation}
see also \cite[$(3.2)$, Prop. 3.2]{Komatsu73}.

\eqref{Avsomega} and \eqref{tildecomparison} yield the important connection between the associated function and the trace function appearing ``naturally'' in the geometric construction.
\end{itemize}

\begin{remark}\label{logconvexmodrem}
\emph{Let $\mathbf{M}=(M_p)_{p\in\NN}\in\RR_{>0}^{\NN}$ such that $\lim_{p\rightarrow+\infty}(M_p)^{1/p}=+\infty$. This approach also implies $\widetilde{\omega}_{\mathbf{M}}\equiv\widetilde{\omega}_{\mathbf{M}^{\on{lc}}}$ and $M_{p_i}=M^{\on{lc}}_{p_i}$ for all $i\in\NN$.}

\emph{In view of \eqref{tildecomparison} and by $M^{\on{lc}}_0=M_0$ one has $\omega_{\mathbf{M}}\equiv\omega_{\mathbf{M}^{\on{lc}}}$ as well and Lemmas \ref{lemma1} and \ref{intformula} can be applied to the sequence $\mathbf{M}^{\on{lc}}$.}
\end{remark}

\subsection{Comparison with Mandelbrojt's approach}\label{comparisonsection}
In \cite[Chapitre I]{mandelbrojtbook} by giving purely geometric arguments and constructions as in Section \ref{normgeomconstrsect} it was stated that
\begin{equation}\label{logconvmin}
\forall\;p\in\NN_{>0}:\;\;\;M^{\on{lc}}_p=\sup_{t\ge 1}\frac{t^p}{\exp(\widetilde{\widetilde{\omega}}_{\mathbf{M}}(t))},
\end{equation}
with
\begin{equation}\label{assoalternative}
\widetilde{\widetilde{\omega}}_{\mathbf{M}}(t):=\sup_{p\in\NN_{>0}}\log\frac{t^p}{M_p},\;\;\;t\ge 1.
\end{equation}
We refer to \cite[p. 17]{mandelbrojtbook} and there the following functions were introduced: $$T:=\exp\circ\widetilde{\widetilde{\omega}}_{\mathbf{M}},\hspace{15pt}A:=\widetilde{\widetilde{\omega}}_{\mathbf{M}}\circ\exp,$$
with $A$ denoting again the \emph{trace function} from \eqref{regsequformula}, \eqref{tracefunction}. However, in \cite{mandelbrojtbook} the function $A$ was only studied on $[0,+\infty)$ which is corresponding to lines having non-negative slopes. Moreover, only indices $p\in\NN_{>0}$ were considered; concerning this issue please see the footnote $(1)$ on \cite[p. 1, Chapitre I]{mandelbrojtbook}.\vspace{6pt}

Of course, the definition of $\widetilde{\widetilde{\omega}}_{\mathbf{M}}$ in \eqref{assoalternative} can be extended to $t\in(0,1]$. But the value $t=0$ is excluded and we see that $\lim_{t\rightarrow 0}\widetilde{\widetilde{\omega}}_{\mathbf{M}}(t)=-\infty$. Thus $\widetilde{\widetilde{\omega}}_{\mathbf{M}}$ also \emph{cannot} be considered on $[0,+\infty)$ as a weight function in the sense of Braun-Meise-Taylor.\vspace{6pt}

However, provided that $\mathbf{M}$ satisfies the regularity assumptions from Lemma \ref{lemma1} we show that for all sufficiently large values the functions $\omega_{\mathbf{M}}$ and $\widetilde{\widetilde{\omega}}_{\mathbf{M}}$ differ in general only by adding a real constant.

\begin{lemma}\label{lemma2}
Let $\mathbf{M}=(M_p)_{p\in\NN}\in\RR_{>0}^{\NN}$ be log-convex and assume that $\lim_{p\rightarrow+\infty}(M_p)^{1/p}=+\infty$. Then we get
\begin{equation}\label{lemma2equ}
\forall\;t\ge\max\{\mu_1,1\}:\;\;\;\omega_{\mathbf{M}}(t)=\log(M_0)+\widetilde{\widetilde{\omega}}_{\mathbf{M}}(t),
\end{equation}
which implies
\begin{equation}\label{lemma2equ1}
\forall\;t\ge\max\{\mu_1,1\}:\;\;\;\widetilde{\omega}_{\mathbf{M}}(t)=\widetilde{\widetilde{\omega}}_{\mathbf{M}}(t).
\end{equation}
\end{lemma}

\demo{Proof}
Let $t\ge 1$ be given such that $\mu_p\le t\le\mu_{p+1}$ for some $p\ge 1$. Then, analogously as in Lemma \ref{lemma1} for any $q>p$ we have $\frac{t^p}{M_p}\ge\frac{t^q}{M_q}$ and, if $p\ge 2$, for any $q<p$ we have $\frac{t^p}{M_p}\ge\frac{t^q}{M_q}$, too. If $p=1$, then according to the definition in \eqref{assoalternative} it is only necessary to consider all $q>p$.

So far we have shown that $\widetilde{\widetilde{\omega}}_{\mathbf{M}}(t)=\log\left(\frac{t^p}{M_p}\right)$ for all $t\ge 1$ satisfying $\mu_p\le t\le\mu_{p+1}$ for some $p\in\NN_{>0}$ and hence \eqref{lemma2equ} follows by taking into account \eqref{lemma1equ}. Finally, \eqref{lemma2equ1} holds by combining \eqref{lemma2equ} and \eqref{tildecomparison}.
\qed\enddemo

\begin{corollary}\label{lemma2cor}
Let $\mathbf{M}\in\hyperlink{LCset}{\mathcal{LC}}$ be given. Then we get $\widetilde{\omega}_{\mathbf{M}}(t)=\omega_{\mathbf{M}}(t)=\widetilde{\widetilde{\omega}}_{\mathbf{M}}(t)$ for all $t\ge\mu_1$.
\end{corollary}

\demo{Proof}
In this case, apart from the basic assumptions on $\mathbf{M}$ in Lemma \ref{lemma2}, we also have $1=M_0\le M_1\Leftrightarrow 1\le\mu_1$ and so \eqref{lemma2equ}, \eqref{lemma2equ1} yield the assertion.
\qed\enddemo

Let now $\mathbf{M}=(M_p)_{p\in\NN}\in\RR_{>0}^{\NN}$ be given such that $\lim_{p\rightarrow+\infty}(M_p)^{1/p}=+\infty$ and set
$$\underline{\underline{M}}_p:=\sup_{t\ge 1}\frac{t^p}{\exp(\widetilde{\widetilde{\omega}}_{\mathbf{M}}(t))}.$$
Then we obtain, compare this with Section \ref{assosequsect}:

\begin{itemize}
\item[$(i)$] First,
$$\sup_{t\ge 1}\frac{t^p}{\exp(\widetilde{\widetilde{\omega}}_{\mathbf{M}}(t))}=\exp\left(\sup_{s\ge 0}\{ps-\widetilde{\widetilde{\omega}}_{\mathbf{M}}(e^s)\}\right)=\exp\left(\varphi^{*}_{\widetilde{\widetilde{\omega}}_{\mathbf{M}}}(p)\right),$$
and so $\underline{\underline{M}}_p$ is log-convex.

\item[$(ii)$] For all $p\in\NN_{>0}$ we get $\underline{\underline{M}}_p\le M_p$ because
$$\underline{\underline{M}}_p=\sup_{t\ge 1}\frac{t^p}{\exp(\widetilde{\widetilde{\omega}}_{\mathbf{M}}(t))}=\sup_{t\ge 1}\frac{t^p}{\sup_{q\in\NN_{>0}}\frac{t^q}{M_q}}\underbrace{\le}_{q=p}M_p.$$

\item[$(iii)$] Finally, if $\mathbf{M}$ is in addition log-convex, then the proof of Lemma \ref{lemma2} yields that
$$\forall\;p\in\NN_{>0},\;\mu_p\ge 1:\;\;\;\underline{\underline{M}}_p=\sup_{t\ge 1}\frac{t^p}{\exp(\widetilde{\widetilde{\omega}}_{\mathbf{M}}(t))}\ge\frac{\mu_p^p}{\exp(\widetilde{\widetilde{\omega}}_{\mathbf{M}}(\mu_p))}=\frac{\mu_p^p}{\mu_p^p/M_p}=M_p;$$
i.e. $\underline{\underline{M}}_p=M_p$ for all $p$ sufficiently large.
\end{itemize}

Summarizing, Corollary \ref{lemma2cor}, Section \ref{assosequsect} and the previous comments yield the following information:

\begin{corollary}\label{doubleunderlinesequcor}
Let $\mathbf{M}\in\hyperlink{LCset}{\mathcal{LC}}$ be given with $1=M_0=M_1$; i.e. $\mu_1=1$. Then
$$\forall\;p\in\NN_{>0}:\;\;\;M^{\on{lc}}_p=M_p=\underline{M}_p=\underline{\underline{M}}_p,$$
and
$$\forall\;p\in\NN_{>0}\;\forall\;\mu_p\le t\le\mu_{p+1}:\;\;\;\omega_{\mathbf{M}}(t)=\widetilde{\omega}_{\mathbf{M}}(t)=\widetilde{\widetilde{\omega}}_{\mathbf{M}}(t)=\log\left(\frac{t^p}{M_p}\right);$$
hence all appearing associated resp. trace functions coincide on $[1,+\infty)$.
\end{corollary}

Moreover, we note:

\begin{itemize}
\item[$(*)$] In Section \ref{lcminorantconstrsect} we assumed $\lim_{p\rightarrow+\infty}(M_p)^{1/p}=+\infty$ which implies, in particular, the existence of some $q_0\in\NN_{>0}$ such that $M_p\ge 1$ for all $p\ge q_0$.

\item[$(*)$] In this case, by changing if necessary the finitely many $M_p$, with $0\le p\le q_0-1$, we can switch to an equivalent sequence $\mathbf{M}^n$ such that $M^n_p\ge 1$ for all $p\in\NN$; i.e. each $S_p=(p,a_p)=(p,\log(M^n_p))$ belongs to the first quadrant. We can set
\begin{equation}\label{normsequdef}
M^n_p:=1,\;\;\;\forall\;0\le p\le q_0-1,\hspace{30pt}M^n_p:=M_p,\;\;\;\forall\;p\ge q_0,
\end{equation}
and both sequences are related by
\begin{equation}\label{strongequ}
\exists\;C\ge 1\;\forall\;p\in\NN:\;\;\;\frac{1}{C}M^n_p\le M_p\le CM^n_p.
\end{equation}

\item[$(*)$] However, even if all points $S_p=(p,a_p)=(p,\log(M_p))$ belong to the first quadrant, then in general it is not sufficient to consider slopes $k\ge 0$ and follow the proof given in \cite{mandelbrojtbook} in order to construct the log-convex minorant: This is due to the missing information of $M_0^{\on{lc}}$; consider e.g. $\mathbf{M}$ given by $M_0=e$, $M_1=1$ and $M_p\ge e$ for all $p\ge 2$.

\item[$(*)$] On the other hand, if each $S_p=(p,a_p)=(p,\log(M_p))$ belongs to the first quadrant or even if $\lim_{p\rightarrow+\infty}(M_p)^{1/p}=+\infty$, then we can replace $\mathbf{M}$ by a sequence $\mathbf{M}^n$ according to \eqref{normsequdef} but such that w.l.o.g. $q_0\ge 2$ and so \eqref{strongequ} still holds. Then one can apply the approach from \cite{mandelbrojtbook} to $\mathbf{M}^n$ and complete the information concerning $p=0$ when setting $(M^n)^{\on{lc}}_0:=(M^n)^{\on{lc}}_1=1$.

\item[$(*)$] \eqref{strongequ} yields by definition
\begin{equation}\label{addrelationweightfct}
\exists\;C\ge 0\;\forall\;t\ge 0:\;\;\;\omega_{\mathbf{M}}(t)-C\le\omega_{\mathbf{M}^n}(t)\le\omega_{\mathbf{M}}(t)+C,
\end{equation}
and via \eqref{tildecomparison} this relation transfers to the corresponding trace functions as well.
\end{itemize}

Gathering all this information, in particular we get:

\begin{proposition}
Let $\mathbf{M}\in\RR_{>0}^{\NN}$ be given such that $\lim_{p\rightarrow+\infty}(M_p)^{1/p}=+\infty$. Then there exists $\mathbf{M}^n\in\RR_{>0}^{\NN}$ which is related to $\mathbf{M}$ by \eqref{strongequ}, such that $(\mathbf{M}^n)^{\on{lc}}\in\hyperlink{LCset}{\mathcal{LC}}$ and such that the corresponding associated functions are related by \eqref{addrelationweightfct}.
\end{proposition}

\section{Non-standard situations for the convex regularization}\label{nonstandardsect}
The aim is to see how the geometric approach described before is changing when $\mathbf{a}$ is violating the basic requirement \eqref{basicregprop} in Section \ref{normgeomconstrsect}. In view of $(b)$ in Section \ref{normgeomconstrsect} we can distinguish between two cases and in both situations formally we could even allow that $M_p=+\infty$ \emph{for infinitely many $p\ge 1$.}

\subsection{Case 1}\label{case1sect}
Let $\mathbf{a}:=(a_p)_{p\in\NN}\in\RR^{\NN}$ be given and assume that
\begin{equation}\label{basicregpropviol}
\liminf_{p\rightarrow+\infty}\frac{a_p}{p}=-\infty;
\end{equation}
i.e. \eqref{basicregprop2} fails. Consequently, by $(b)$ in Section \ref{normgeomconstrsect} also $\lim_{p\rightarrow+\infty}\frac{a_p}{p}=+\infty$ is violated. When using the notation from this Section we see that in this case, since $\frac{a_p}{p}$ is the slope of the straight line connecting the points $(0,0)$ and $S_p$, we cannot find a straight line $D_k$ with $k>-\infty$ such that $D_k$ is passing through the point $S_0$ and such that all $S_p$, $p\ge 1$, are lying strictly above $D_k$.\vspace{6pt}

Thus in the geometric approach we can only consider the $y$-axis as the limiting case which does have ``slope'' $-\infty$. In order to give the regularized sequence $\widetilde{S}=(\widetilde{S}_p)_{p\in\NN}$ a meaning, we have to put
$$\widetilde{a}_0=a_0,\hspace{30pt}\widetilde{a}_p=-\infty,\;\;\;\forall\;p\in\NN_{>0}.$$
In this case the set $\{D_k: k\in\RR\}$ is empty and formally one can only consider $D_{-\infty}$ which corresponds to the set $\{(0,y): y\in\RR\}$. $S_0$ is the only point belonging to $\mathcal{S}$ and lying on this line. Thus, the trace function in \eqref{tracefunction} takes the form $A(-\infty)=-d_k=-a_0$ and it is not defined for $k\in\RR$. Consequently, $\widetilde{a}_0=a_0$ but $\widetilde{a}_p<a_p$ for all $p\ge 1$ and hence the sequences $\mathbf{a}$, $\widetilde{\mathbf{a}}$ differ at all indices $p\in\NN_{>0}$. Indeed, we have $\widetilde{\mathbf{a}}\notin\RR^{\NN}$.\vspace{6pt}

Moreover, when $\mathbf{M}=(M_p)_{p\in\NN}\in\RR_{>0}^{\NN}$ is given, then via \eqref{logconvtrafo} condition \eqref{basicregpropviol} precisely means
$$\mathbf{M}_{\iota}=\liminf_{p\rightarrow+\infty}(M_p)^{1/p}=0;$$
i.e. \eqref{liminfcond} is violated. In this case, by \eqref{logconvtrafo1}, the sequence $\mathbf{M}^{\on{lc}}$ is given by
$$M^{\on{lc}}_0=M_0,\hspace{30pt}\forall\;p\in\NN_{>0}:\;\;\;M^{\on{lc}}_p=0.$$
Using and combining this information, formally we see that the formulas in Section \ref{lcminorantconstrsect} involving the trace function $A$, and so $\widetilde{\omega}_{\mathbf{M}}$, $\omega_{\mathbf{M}}$, for computing $\mathbf{M}^{\on{lc}}$ fail because here we have $\widetilde{\omega}_{\mathbf{M}}(0)=A(-\infty)=-a_0=-\log(M_0)$ and $A$ is not defined on $\RR$. However, recall $(ii)$ in Lemma \ref{assofunctsectlemma}: We have $\omega_{\mathbf{M}}(0)=0$ and $\omega_{\mathbf{M}}(s)=+\infty$ for all $s>0$ and so \eqref{tildecomparison} formally only makes sense for $s=0$. But when ``extending'' this identity now to $(0,+\infty)$ by setting $\widetilde{\omega}_{\mathbf{M}}(s)=+\infty$ for all $s\in(0,+\infty)$, using this and the conventions $0^0:=1$ and $\frac{1}{+\infty}=0$ too, then the formula $M^{\on{lc}}_p=\sup_{t\ge 0}\frac{t^p}{\exp(\widetilde{\omega}_{\mathbf{M}}(t))}$ makes sense for all $p\in\NN$. Finally, via \eqref{Avsomega} we can set $A(k):=+\infty$ for all $k\in\RR$.

\subsection{Case 2}\label{case2sect}
Let $\mathbf{a}:=(a_p)_{p\in\NN}\in\RR^{\NN}$ be given and assume that
\begin{equation}\label{basicregpropvio2}
-\infty<\mathbf{a}_{\iota}:=\liminf_{p\rightarrow+\infty}\frac{a_p}{p}<+\infty.
\end{equation}
We use again the notation from Section \ref{normgeomconstrsect} and start with the geometric construction as described there since we can find straight lines $D_k$ with $k>-\infty$ and such all $S_p$, $p\ge 1$, are lying strictly above $D_k$ for all $k\in\RR$ sufficiently small. Then we proceed with the construction described in Section \ref{normgeomconstrsect}. Since $\frac{a_p}{p}$ is the slope of the straight line connecting $(0,0)$ and $S_p$ and since
\begin{equation}\label{aiotadef}
\mathbf{a}_{\iota}<+\infty,
\end{equation}
we have that \eqref{regsequformula} turns into
\begin{equation}\label{regsequformulacase2}
\forall\;p\in\NN:\;\;\;\widetilde{a}_p=\sup_{k\in(-\infty,\mathbf{a}_{\iota})}\{kp+d_k\}=\sup_{k\in(-\infty,\mathbf{a}_{\iota})}\{kp-A(k)\},
\end{equation}
where again $A(k)=-d_k$ is denoting the \emph{trace function;} see \eqref{tracefunction}. Note that also here $A$ is continuous and non-decreasing, that $\widetilde{a}_0=\sup_{k\in(-\infty,\mathbf{a}_{\iota})}\{d_k\}=a_0$ and we put again $A(-\infty):=-a_0$. By the conventions $0\cdot(-\infty):=0$ and $p\cdot(-\infty):=-\infty$ for any $p\in\NN_{>0}$ we see that in \eqref{regsequformulacase2} one can consider $k\in\{-\infty\}\cup(-\infty,\mathbf{a}_{\iota})$.

\eqref{dkequ} transfers to $d_k=\inf_{p\in\NN}\{a_p-pk\}$ for all $k\in(-\infty,\mathbf{a}_{\iota})$ and consequently \eqref{trace} takes the following form:
\begin{equation*}\label{tracecase2}
\forall\;k\in(-\infty,\mathbf{a}_{\iota}):\;\;\;A(k)=-d_k=-\inf_{p\in\NN}\{a_p-pk\}=\sup_{p\in\NN}\{pk-a_p\}.
\end{equation*}
$A$ is now naturally defined on the open interval $(-\infty,\mathbf{a}_{\iota})$.

Given $\mathbf{M}\in\RR_{>0}^{\NN}$ we apply this to the sequence $a_p:=\log(M_p)$ and so \eqref{basicregpropvio2} precisely means that $$0<\mathbf{M}_{\iota}=\liminf_{p\rightarrow+\infty}(M_p)^{1/p}=\exp(\mathbf{a}_{\iota})<+\infty.$$
Consequently, we get for any $k\in(-\infty,\mathbf{a}_{\iota})$ that
$$A(k)=\sup_{p\in\NN}\{pk-a_p\}=\sup_{p\in\NN}\{pk-\log(M_p)\}=\sup_{p\in\NN}\log\left(\frac{e^{pk}}{M_p}\right)=\widetilde{\omega}_{\mathbf{M}}(e^k),$$
thus \eqref{regsequformulacase2}, recall also \eqref{logconvtrafo} and \eqref{logconvtrafo1}, yields
\begin{align*}
\forall\;p\in\NN:\;\;\;M^{\on{lc}}_p&=\exp(\widetilde{a}_p)=\exp\left(\sup_{k\in(-\infty,\mathbf{a}_{\iota})}\{kp-A(k)\}\right)=\sup_{k\in(-\infty,\mathbf{a}_{\iota})}\frac{e^{kp}}{\exp(A(k))}
\\&
=\sup_{k\in(-\infty,\mathbf{a}_{\iota})}\frac{e^{kp}}{\exp(\widetilde{\omega}_{\mathbf{M}}(e^k))}=\sup_{s\in(0,\exp(\mathbf{a}_{\iota}))}\frac{s^p}{\exp(\widetilde{\omega}_{\mathbf{M}}(s))}=\sup_{s\in(0,\mathbf{M}_{\iota})}\frac{s^p}{\exp(\widetilde{\omega}_{\mathbf{M}}(s))}.
\end{align*}
We summarize:

\begin{itemize}
\item[$(a)$] $M^{\on{lc}}_0=M_0$ holds since $$\lim_{s\rightarrow0}\widetilde{\omega}_{\mathbf{M}}(s)=\lim_{k\rightarrow-\infty}A(k)=\lim_{k\rightarrow-\infty}-d_k=-a_0=-\log(M_0).$$ Moreover, in the supremum above we can consider all $s\in[0,\mathbf{M}_{\iota})$. Using the correspondence with the trace function $A$ we see that \eqref{tildecomparison} is valid for all $s\in[0,\mathbf{M}_{\iota}]$ with considering the limit at the right-end point of this interval. Note that, analogously to the previous Section \ref{case1sect}, by taking into account $(ii)$ in Lemma \ref{assofunctsectlemma}, formally we can put $\widetilde{\omega}_{\mathbf{M}}(s)=+\infty$ for all $s>\mathbf{M}_{\iota}$, i.e. set $A(k)=+\infty$ for all $k\in(\mathbf{a}_{\iota},+\infty)$, and use the convention $\frac{1}{+\infty}=0$ in order to consider $\sup_{s\ge 0}$ instead of $\sup_{s\in(0,\mathbf{M}_{\iota})}$ in the above representation for $M^{\on{lc}}_p$.

\item[$(b)$] We have the estimate
\begin{equation}\label{lcequ}
\forall\;p\in\NN:\;\;\;M^{\on{lc}}_p=\sup_{s\in(0,\mathbf{M}_{\iota})}\frac{s^p}{\exp(\widetilde{\omega}_{\mathbf{M}}(s))}\le\frac{\mathbf{M}_{\iota}^p}{\exp(-a_0)}=M_0\mathbf{M}_{\iota}^p,
\end{equation}
since $\widetilde{\omega}_{\mathbf{M}}$ is non-decreasing and for $p=0$ we even have equality.

Thus, up to a constant, $\mathbf{M}^{\on{lc}}$ is lying below a sequence which is increasing geometrically if $\mathbf{M}_{\iota}>1$ and decreasing geometrically if $0<\mathbf{M}_{\iota}<1$. If $\mathbf{M}_{\iota}=1$, then $\mathbf{M}^{\on{lc}}$ is uniformly bounded by $M_0$.

\item[$(c)$] \eqref{basicregpropvio2} is not excluding $\mathbf{M}_{\sigma}:=\limsup_{p\rightarrow+\infty}(M_p/M_0)^{1/p}=\limsup_{p\rightarrow+\infty}(M_p)^{1/p}=+\infty$ and in this case for infinitely many indices $p$ the difference between $M^{\on{lc}}_p$ and $M_p$ can be seen to be very large: More precisely, by \eqref{lcequ} one has
    \begin{equation}\label{largedifferencequ}
    \sup_{p\in\NN_{>0}}\left(\frac{M_p}{M^{\on{lc}}_p}\right)^{1/p}\ge\sup_{p\in\NN_{>0}}\frac{(M_p/M_0)^{1/p}}{\mathbf{M}_{\iota}}=+\infty,
    \end{equation}
thus $\mathbf{M}^{\on{lc}}$ and $\mathbf{M}$ are not equivalent.
\end{itemize}

Comment $(c)$ yields the following consequence.

\begin{corollary}
Let $\mathbf{M}\in\RR_{>0}^{\NN}$ such that
$$0<\mathbf{M}_{\iota}<\mathbf{M}_{\sigma}=+\infty.$$
Then there does not exist a log-convex sequence $\mathbf{L}\in\RR_{>0}^{\NN}$ such that $\mathbf{M}$ and $\mathbf{L}$ are equivalent.
\end{corollary}

This result is in contrast to the case where $$0<\mathbf{M}_{\iota}<\mathbf{M}_{\sigma}<+\infty,$$
because here $C^{p+1}\le M_p\le D^{p+1}$ for some constants $0<C<1<D$ and all $p\in\NN$ and so $\mathbf{M}$ is equivalent to the constant sequence $\mathbf{1}:=(1)_{p\in\NN}$.

\demo{Proof}
Let $\mathbf{L}$ be a log-convex sequence such that $\mathbf{M}$ and $\mathbf{L}$ are equivalent. So $\frac{1}{C^{p+1}}L_p\le M_p$ for some $C\ge 1$ and all $p\in\NN$ is satisfied. Obviously, $\mathbf{L}^{1/C}:=\left(\frac{1}{C^{p+1}}L_p\right)_{p\in\NN}$ is also log-convex and since $\mathbf{L}^{1/C}\le\mathbf{M}$ this implies by definition of $\mathbf{M}^{\on{lc}}$ that $\mathbf{L}^{1/C}\le\mathbf{M}^{\on{lc}}$, too. However, by assumption clearly $\mathbf{L}^{1/C}$ is equivalent to $\mathbf{M}$ as well, a contradiction to \eqref{largedifferencequ}.
\qed\enddemo

\begin{proposition}
Let $\mathbf{M}\in\RR_{>0}^{\NN}$ such that $0<\mathbf{M}_{\iota}<+\infty$. Then, with $\mu^{\on{lc}}_p:=\frac{M^{\on{lc}}_p}{M^{\on{lc}}_{p-1}}$ we get that $$\lim_{p\rightarrow+\infty}\mu^{\on{lc}}_p=\lim_{p\rightarrow+\infty}(M^{\on{lc}}_p)^{1/p}=\mathbf{M}_{\iota}$$
is valid and hence Lemma \ref{lemma1exotic} can be applied to $\mathbf{M}^{\on{lc}}$ and to $C:=\mathbf{M}_{\iota}$.
\end{proposition}

\demo{Proof}
In view of Lemma \ref{preliminarysectlemma} and the comments in Section \ref{preliminarysection} applied to $\mathbf{M}^{\on{lc}}$ we have that both limits $\lim_{p\rightarrow+\infty}\mu^{\on{lc}}_p$ and $\lim_{p\rightarrow+\infty}(M^{\on{lc}}_p)^{1/p}$ exist and coincide.

Let us show that $\lim_{p\rightarrow+\infty}(M^{\on{lc}}_p)^{1/p}=\mathbf{M}_{\iota}=\liminf_{p\rightarrow+\infty}(M_p)^{1/p}$:

Since $(M^{\on{lc}}_p)^{1/p}\le(M_p)^{1/p}$ for all $p\ge 1$ we clearly get $\lim_{p\rightarrow+\infty}(M^{\on{lc}}_p)^{1/p}\le\liminf_{p\rightarrow+\infty}(M_p)^{1/p}$. If $\lim_{p\rightarrow+\infty}(M^{\on{lc}}_p)^{1/p}<\mathbf{M}_{\iota}$, then there exist $C>C_1>0$ such that $M^{\on{lc}}_p<C_1^p<C^p<M_p$ for all $p$ sufficiently large and choosing $D\ge 1$ large enough we have $D^{-1}C^p\le M_p$ for all $p\in\NN$. The auxiliary sequence $\mathbf{L}:=(D^{-1}C^p)_{p\in\NN}$ is clearly log-convex, satisfies $\mathbf{L}\le\mathbf{M}$ and so $\mathbf{L}\le\mathbf{M}^{\on{lc}}$ holds. But this yields $D^{-1}C^p\le M^{\on{lc}}_p<C_1^p$ for all $p$ sufficiently large, a contradiction as $p\rightarrow+\infty$.
\qed\enddemo

We present now examples which illustrate the fact that in general we do not have $\widetilde{a}_p=a_p$ resp. $M^{\on{lc}}_p=M_p$ for infinitely many $p$; i.e. this is in contrast to the approach studied in Section \ref{lcminorsect}.

\begin{example}\label{nonstandardsectex}
\emph{We consider ``almost affine linear'' sequences.}
\begin{itemize}
\item[$(i)$] \emph{Let $\mathbf{a}=(a_p)_{p\in\NN}$ be given as follows:}
\begin{equation}\label{nonstandardsectexequ}
a_0:=0,\hspace{30pt}a_1:=-1,\hspace{30pt}a_p:=cp,\;\;\;p\ge 2,\;c\ge 0.
\end{equation}
\emph{Then $\widetilde{a}_0=a_0$, $\widetilde{a}_1=a_1$, but $\widetilde{a}_p<a_p$ for all $p\ge 2$. Indeed, for all $p\ge 2$ the points $S_p$ are projected vertically down onto the line $D_{c}^{-}:=\lim_{\varepsilon\rightarrow 0}D_{c-\varepsilon}$ with $D_{c-\varepsilon}$ being the straight line passing through $S_1=(1,-1)$ and having slope $c-\varepsilon$.}

\item[$(ii)$] \emph{On the other hand, if $a_1:=c$ instead and all other terms remain unchanged, then $\widetilde{\mathbf{a}}\equiv\mathbf{a}$ and so the convex minorant coincides with $\mathbf{a}$. Note that the points $S_p=(p,a_p)$ are lying on the straight line having slope $c$ and passing through $S_0=(0,0)$ and so $\mathbf{a}$ is clearly convex. However, the difference compared with Section \ref{normgeomconstrsect} is that ``during the geometric construction'' we do not get infinitely many points where the regularized sequence coincides with the given one and the equality, for infinitely many indices, is finally obtained in the limit: Each $S_p$ is projected vertically down onto the line $D_{c}^{-}:=\lim_{\varepsilon\rightarrow 0}D_{c-\varepsilon}$ with $D_{c-\varepsilon}$ denoting the straight line passing through $(0,0)$ and having slope $c-\varepsilon$. But $S_0=(0,0)$ is the only point belonging to the set $\mathcal{S}$ and lying on $D_{c-\varepsilon}$, $\varepsilon>0$ small but arbitrary.}
\end{itemize}
\emph{In both cases \eqref{basicregpropvio2} is valid with the same value since $\lim_{p\rightarrow+\infty}\frac{a_p}{p}=c=\mathbf{a}_{\iota}$. Via \eqref{Mvsarelation} we see that $(i)$ corresponds to $M_0=1$, $M_1=e^{-1}$, and $M_p=e^{cp}$, $p\ge 2$, resp. $M_1=e^c$ instead in $(ii)$. In both cases $\mathbf{M}_{\iota}=e^c$, but in $(i)$ we have $\mathbf{M}_{\inf}=e^{-1}<\mathbf{M}_{\iota}$ whereas $\mathbf{M}_{\inf}=e^{c}=\mathbf{M}_{\iota}$ in $(ii)$.}

\emph{Using this we show that in general $\mathbf{M}_{\inf}$ cannot be replaced by $\mathbf{M}_{\iota}$ in $(i)$ in Lemma \ref{assofunctsectlemma}: The first example yields $A(k)=0$ for $k\in(-\infty,-1]$ but then $A(k)\rightarrow c+1$ as $k\rightarrow c$ resp. $\widetilde{\omega}_{\mathbf{M}}(s)\rightarrow c+1$ as $s\rightarrow\mathbf{M}_{\iota}$ because $(0,-1-c)$ is the intersecting point of the $y$-axis with the straight line having slope $c$ and passing through $S_1=(1,-1)$. More precisely, $A(k)>0$ for all $k\in(-1,\mathbf{a}_{\iota})$ resp. $\widetilde{\omega}_{\mathbf{M}}(s)>0$ for all $s\in(e^{-1},\mathbf{M}_{\iota})$. By taking into account \eqref{tildecomparison} we get $\omega_{\mathbf{M}}(s)\rightarrow c+1$ as $s\rightarrow\mathbf{M}_{\iota}$ and $\omega_{\mathbf{M}}(s)>0$ for all $\mathbf{M}_{\inf}<s\le\mathbf{M}_{\iota}$.}
\end{example}

We close by showing that in the situation studied in this Section it can also appear that one has equality between $\widetilde{\mathbf{a}}$ and $\mathbf{a}$ for infinitely many indices.

\begin{example}\label{nonstandardsectex1}
\emph{Let $a_0:=0$ and each $a_p$, $p\ge 1$, is assumed to lie on the straight line passing through $S_0=(0,0)$ and having slope $c_p$; i.e. $a_p=c_pp$. Let $c_1:=1$ and fix $c>1$. We claim to define iteratively a strictly increasing sequence $(c_p)_{p\ge 1}$ such that}
\begin{equation}\label{sloperelation}
(p+1)c_{p+1}-pc_p<c,\;\;\;p\ge 1,\hspace{30pt}\frac{2(p+1)}{p+2}c_{p+1}-\frac{p}{p+2}c_p<c_{p+2},\;\;\;p\ge 1.
\end{equation}
\emph{The first relation gives that the slope of the line connecting $S_p=(p,a_p)$ and $S_{p+1}=(p+1,a_{p+1})$ is strictly smaller than the value $c$ for any $p\ge 1$. Note that the case for $S_0$ and $S_1$ is clear since this slope is given by $c_1=1<c$. The second relation in \eqref{sloperelation} precisely means that the slope of the line connecting $S_p$ and $S_{p+1}$ is strictly smaller than the one connecting $S_{p+1}$ and $S_{p+2}$ for all $p\ge 1$. For $p=0$ we have to compare the slope $c_1$ with $2c_2-c_1$ and $c_1<2c_2-c_1\Leftrightarrow c_1<c_2$ since the sequence is strictly increasing as verified below.}\vspace{6pt}

\emph{Consequently, when applying the geometric approach from before, for which it suffices to consider non-negative slopes, then in fact $\widetilde{a}_p=a_p$ for all $p\in\NN$ but $\limsup_{p\rightarrow+\infty}\frac{a_p}{p}\le c<+\infty$.}\vspace{6pt}

\emph{We introduce $(c_p)_{p\ge 1}$ inductively by
$$c_{p+1}:=\frac{c}{2(p+1)}+\frac{2p+1}{2(p+1)}c_p,\;\;\;p\ge 1;$$
i.e. $c_{p+1}$ denotes the middle point of the interval $[c_p,\frac{c}{p+1}+\frac{p}{p+1}c_p=:d_p]$.}

\emph{First $c_p<c$ for all $p\ge 1$ since this is clear for $p=1$ and $c_{p+1}=\frac{c}{2(p+1)}+\frac{2p+1}{2(p+1)}c_p<c\Leftrightarrow(2p+1)c_p<(2p+1)c\Leftrightarrow c_p<c$ holds then by induction. Similarly, one can prove that $d_p<c\Leftrightarrow c_p<c$ and, moreover, $c_p<d_p\Leftrightarrow c_p<c$ holds and hence inductively the choice for $c_{p+1}$ makes sense for all $p\ge 1$.}

\emph{The sequence $(c_p)_p$ is strictly increasing since $c_p<c_{p+1}\Leftrightarrow c_p<\frac{c}{2(p+1)}+\frac{2p+1}{2(p+1)}c_p\Leftrightarrow c_p<c$.}

\emph{By definition, the first requirement in \eqref{sloperelation} is equivalent to $c_p<c$ and hence valid. Let us check the second one: We have}
\begin{align*}
\frac{2(p+1)}{p+2}c_{p+1}-\frac{p}{p+2}c_p&=\frac{2(p+1)}{p+2}\frac{c}{2(p+1)}+\frac{2(p+1)}{p+2}\frac{2p+1}{2(p+1)}c_p-\frac{p}{p+2}c_p
\\&
=\frac{c}{p+2}+\frac{2p+1}{p+2}c_p-\frac{p}{p+2}c_p=\frac{4pc+4c}{4(p+1)(p+2)}+\frac{4(p+1)(p+1)}{4(p+1)(p+2)}c_p
\end{align*}
\emph{and}
\begin{align*}
c_{p+2}&=\frac{c}{2(p+2)}+\frac{2p+3}{2(p+2)}c_{p+1}=\frac{c}{2(p+2)}+\frac{2p+3}{2(p+2)}\frac{c}{2(p+1)}+\frac{2p+3}{2(p+2)}\frac{2p+1}{2(p+1)}c_p
\\&
=\frac{4cp+5c+(2p+3)(2p+1)c_p}{4(p+2)(p+1)}.
\end{align*}
\emph{Thus the required estimate is equivalent to
$4pc+5c+(2p+3)(2p+1)c_p>4pc+4c+4(p+1)^2c_p\Leftrightarrow c+(4p^2+8p+3)c_p>(4p^2+8p+4)c_p\Leftrightarrow c>c_p$ and so this holds inductively as shown above for all $p\ge 1$, too.}
\end{example}

\section{Regularizing sequences by a regularizing function}\label{generalsection}
We discuss now the more general approach of regularizing a sequence $\mathbf{a}$ w.r.t. a regularizing function $\phi$. Let $\mathbf{a}:=(a_p)_{p\in\NN}\in\RR^{\NN}$ be given and set again $S_p:=(p,a_p)$, $p\in\NN$, and $\mathcal{S}:=\{S_p: p\in\NN\}$ as in Section \ref{normgeomconstrsect}. For the whole procedure described in this Section we can formally allow again that $a_p=+\infty$ \emph{for finitely many indices $p\ge 1$;} when $\lim_{t\rightarrow T}\phi(T)=+\infty$ for some $T\in\RR$, then also \emph{for infinitely many $p\ge 1$} or \emph{even for all $p\ge 1$ sufficiently large.}

\subsection{Regularizing function}\label{regufunctionsect}
\begin{definition}\label{regufctdef}
A function $\phi:\RR\rightarrow\RR_{\ge 0}\cup\{+\infty\}$ is called a \emph{regularizing function} if:
\begin{itemize}
\item[$(I)$] $\phi$ is non-decreasing,

\item[$(II)$] $\lim_{t\rightarrow-\infty}\phi(t)=0$,

\item[$(III)$] $\lim_{t\rightarrow T}\phi(t)=+\infty$ for $T\in\RR$ or $T=+\infty$,

\item[$(IV)$] $\phi$ is continuous on $(-\infty,T)$ for $T\in\RR$ or $T=+\infty$, with $T$ denoting the value from $(III)$.
\end{itemize}
We write $\phi(T)=+\infty$ for some $T\in\RR$ if $\lim_{t\rightarrow T}\phi(T)=+\infty$ and $\phi\neq+\infty$ if $\phi(t)<+\infty$ for all $t\in\RR$; i.e. if $T=+\infty$.
\end{definition}

Note: If $\phi(T)=+\infty$ for some $T\in\RR$, then we set $\phi(t)=+\infty$ for all $t\ge T$. We discuss now the geometric process of regularizing the given sequence $\mathbf{a}\in\RR^{\NN}$ or the set $\mathcal{S}$ w.r.t. a regularizing function $\phi$. The regularized sequence is then denoted by $\mathbf{a}^{\phi}:=(a^{\phi}_p)_{p\in\NN}$; similarly write $S^{\phi}_p$ and $\mathcal{S}^{\phi}$. We work geometrically in $\RR^2$ with coordinates $(x,y)$. A special case, formally not within the scope of the basic properties $(I)-(IV)$, is $\phi=+\infty$; i.e. $\phi(t)=+\infty$ for all $t\in\RR$. The forthcoming construction gives resp. reproves in this special situation the approach from Sections \ref{lcminorsect} and \ref{nonstandardsect}.

\subsection{Geometric notions}\label{geomnotsect}
The main idea in order to proceed is to consider the following purely geometric construction depending on $\phi$ and $\mathbf{a}$:\vspace{6pt}

\begin{itemize}
\item[$(1)$] First, we introduce the set
$$B_t:=\{(x,y)\in\RR^2: 0\le x\le\phi(t)\},\;\;\;t\in\RR;$$
i.e. $B_t$ is the stripe on the right-half plane in $\RR^2$ with width $\phi(t)$.

\item[$(2)$] Second, for all $t\in\RR$ we define $D_t$ to be the intersection of $B_t$ with the straight line having slope $t$ and satisfying the following two properties:
    \begin{itemize}
    \item[$(2.1)$] There exists at least one element in the set $\mathcal{S}\cap B_t$ such that this point also belongs to $D_t$ and

    \item[$(2.2)$] such that none of the points $S_p\in B_t$ is lying below $D_t$.
    \end{itemize}

\emph{In other words:} When $t\in\RR$ is given and fixed, then we move the straight line with slope $t$ parallel ``from below'' until it touches at least one of the points $S_p\in B_t$. Then the intersection of this line with $B_t$ is denoted by $D_t$.

Occasionally, we consider the ``limiting case'' $B_{-\infty}$ resp. $D_{-\infty}$ and both sets correspond to the $y$-axis.
\end{itemize}

First let us comment on the following fact: The previously defined setting is equivalent to the even more general situation when we consider in the stripe $B_t$ the segment $D_{\psi(t)}$ having slope $\psi(t)$ with $\psi:\RR\rightarrow\RR$ being a strictly increasing and continuous function such that $\lim_{t\rightarrow-\infty}\psi(t)=-\infty$, $\lim_{t\rightarrow+\infty}\psi(t)=+\infty$. This follows because in this more general case we consider the reparametrization $\phi\circ\psi^{-1}$, which is clearly again a regularizing function, and apply the framework in this Section to $\phi\circ\psi^{-1}$. Note that such a reparametrization does not change the crucial different behavior of the regularizing function: $\phi=+\infty$, $\phi\neq+\infty$, and $\phi(T)=+\infty$ for some $T\in\RR$.\vspace{6pt}

Next let us see that the mapping $t\mapsto D_t$ is well-defined in the sense that for each given $t\in\RR$ we can find at least one point $S_p\in D_t$.

\begin{itemize}
\item[$(*)$] If there exists some $T\in\RR$ such that $\phi(t)=0$ for all $t\in(-\infty,T]$, then $B_t$ coincides with the $y$-axis for all $t\in(-\infty,T]$. In this case for all such $t$ the segment $D_t$ consists only of its intersection point with the $y$-axis; i.e. we have $D_t=S_0$. Note that $S_0=(0,a_0)$ and $a_0>-\infty$.

\item[$(*)$] For any regularizing function $\phi$, since $\phi\ge 0$, at least $S_0$ belongs to each $B_t$ and hence is lying on $D_t$.

\item[$(*)$] If $\phi\neq +\infty$, then for all $t\in\RR$ with $\phi(t)<1$ the point $S_0$ is indeed the only one and whenever $\phi(t)\ge 1$ then (at least) a second point $S_p$, $p\ge 1$, belongs to $B_t$. This comment should be compared with the ``normalization condition'' $\phi(0)\ge 1$ used in \cite{mandelbrojtbook}.

    For this fact it is crucial to assume naturally that $a_p>-\infty$ for all $p\in\NN$ and which is ensured since $\mathbf{a}\in\RR^{\NN}$: Otherwise, if there exists $p_0\in\NN$ with $a_{p_0}=-\infty$, then for all $t\in\RR$ such that $\phi(t)\ge p_0$ the geometric idea is not working since we would have to put $D_t: y=-\infty$ in order to give this case a meaning. However, if ``only'' \eqref{basicregprop1} is violated, i.e. if
    \begin{equation}\label{basicregprop1viol}
    \liminf_{p\rightarrow+\infty}a_p=-\infty,
    \end{equation}
    then as $t\rightarrow+\infty$ the lines $D_t$ are passing through points whose $y$-coordinate is approaching the value $-\infty$. So even this case can be treated for such regularizing functions and the case if \eqref{basicregprop2} is violated is similar.

\item[$(*)$] If $\phi=+\infty$, then each $B_t$ is the whole right half-plane/space and so $D_t$ is a half-ray and we are allowed to use half-rays with arbitrary slope $t\in\RR$. Thus this case corresponds to the situation treated in Sections \ref{lcminorsect} and \ref{nonstandardsect}. Hence, for $t\mapsto D_t$ being well-defined one should avoid \eqref{basicregpropviol} since in this case for any $t\in\RR$ at least one point will lie strictly below $D_t$; see Section \ref{case1sect} for this degenerate case.

\item[$(*)$] If $\phi(T)=+\infty$ for some $T\in\RR$, then as $t\rightarrow T$ the set $B_t$ is tending to the right-half plane and so $D_t$ to a half-ray having slope $T$. Here it is also clear that $t\mapsto D_t$ is well-defined for any $\mathbf{a}\in\RR^{\NN}$ and, as we are going to see, this situation allows to treat formally also sequences $\mathbf{a}$ such that $a_p=+\infty$ for all $p\in\NN_{>0}$ sufficiently large. Similarly as for the case $\phi\neq+\infty$, even if \eqref{basicregprop1} resp. \eqref{basicregprop2} is violated, the construction makes sense.
\end{itemize}

\subsection{Principal indices}\label{principalsect}
Now consider the subsequence $(S_{p_i})_{i\in\NN}$ of points with the following property: There exists some $t\in\RR$ such that $S_{p_i}$ is lying on $D_t$. The corresponding sequence of indices $(p_i)_{i\in\NN}$ is called the sequence of \emph{principal indices} and $S_{p_i}$ is called a \emph{principal point.}

\begin{itemize}
\item[$(i)$] $S_{p_0}=S_0$ holds for any regularizing function $\phi$.

\item[$(ii)$] If $\phi\neq +\infty$, then $(S_{p_i})_{i\in\NN}$ consists of infinitely many points: If there exists an index $i_0\in\NN$ such that $S_p\notin D_t$ for any $t\in\RR$ and for all $p>p_{i_0}$, then $\phi(t)\rightarrow+\infty$ as $t\rightarrow+\infty$, and hence we are allowed to consider arbitrary large slopes. Consequently, we have $a_p=+\infty$ for all $p>p_{i_0}$. But this is a contradiction to the basic assumption that $a_p=+\infty$ for possibly at most finitely many $p\ge 1$.

\item[$(iii)$] If $\phi=+\infty$, then in order to guarantee the existence of infinitely many principal indices we have to assume \eqref{basicregprop}; i.e. the standard approach treated in Section \ref{lcminorsect}. If \eqref{basicregprop} is violated, then in general one only has finitely many indices $p_i$; see Section \ref{nonstandardsect} and, in particular, Example \ref{nonstandardsectex}: In $(i)$ there we can find two indices $p_0=0$ and $p_1=1$, whereas only $p_0=0$ in $(ii)$. In the situation treated in Section \ref{case1sect} there exists only one principal index, namely $p_0=0$. Concerning Section \ref{case2sect} we point out that the behavior of the sequence $(a_p/p)_{p\in\NN_{>0}}$ near $\mathbf{a}_{\iota}$ is crucial, see also Example \ref{nonstandardsectex1}: For sure there exist finitely many principal indices if $\frac{a_p}{p}\ge\mathbf{a}_{\iota}$ for all but possibly finitely many values $p\ge 1$.

\item[$(iv)$] If $\phi(T)=+\infty$ for some $T\in\RR$, then it is also not clear that in general infinitely many principal indices do exist even if $\sup_{p\in\NN}a_p<+\infty$: Assume that $p_{i_0}$ is a principal index, and there exists at least one, namely $p_0=0$, but all points $S_p$ with $p\ge i_0+1$ are lying on or even strictly above the straight line with slope $T$ and which is passing through $S_{p_{i_0}}$. This line is given by
    \begin{equation}\label{lineT}
    \ell_T: s\mapsto sT+d_T=sT+a_{p_{i_0}}-Tp_{i_0}=T(s-p_{i_0})+a_{p_{i_0}},
    \end{equation}
    and so the situation described before holds if
    $$\forall\;p\ge p_{i_0}+1:\;\;\;a_p\ge T(p-p_{i_0})+a_{p_{i_0}}.$$
    As $t\rightarrow T$ we are only allowed to approach $\ell_T$, formally $\ell_T=D_{T}^{-}:=\lim_{\varepsilon\rightarrow 0}D_{T-\varepsilon}$, and hence all $S_p$, $p\ge i_0+1$, will lie strictly above all segments $D_t$ under consideration and $p_{i_0}$ is the ``last principal index''.\vspace{6pt}

    If there exist only finitely many principal indices then we say that eventually $\mathbf{a}$ lies on lines having at least slope $T$.

\begin{itemize}
   \item[$(*)$] On the one hand, if $\mathbf{a}$ has this property, then $\frac{a_p}{p}\ge T(1-\frac{p_{i_0}}{p})+\frac{a_{p_{i_0}}}{p}$ for all $p\ge p_{i_0}+1$ with $p_{i_0}$ denoting the last principal index, and hence
    \begin{equation}\label{liminfT0}
    \liminf_{p\rightarrow+\infty}\frac{a_p}{p}\ge T.
    \end{equation}
   Thus,
    \begin{equation}\label{liminfT}
    \liminf_{p\rightarrow+\infty}\frac{a_p}{p}<T,
    \end{equation}
    is sufficient to ensure the existence of infinitely many principal indices $p_i$.

   \item[$(*)$] On the other hand, if in \eqref{liminfT0} we have a strict inequality, then we can find $\epsilon>0$ and $p_{\epsilon}\in\NN$ such that $a_p\ge Tp+\epsilon p$ for all $p\ge p_{\epsilon}$. So eventually each $S_p$ is lying strictly above the straight line with slope $T$ and passing through the origin. In this case there exist finitely many principal indices.
\end{itemize}

For sure there exist only finitely many principal indices if $a_p=+\infty$ for all $p\in\NN_{>0}$ sufficiently large; a formal situation which can be considered only for this case. Moreover, in this setting we can assume that for infinitely many $p\ge 1$ we have $a_p=+\infty$.\vspace{6pt}

Summarizing, $\phi(T)=+\infty$ for some $T\in\RR$ is very similar to $\phi=+\infty$ in Section \ref{case2sect}: In the first case the admissible slopes are uniformly bounded by $T$ and in the latter case by the value $\mathbf{a}_{\iota}$. However, in Section \ref{case2sect} formally we are allowed to consider $a_p=+\infty$ for infinitely many $p\ge 1$ but \emph{not for all $p\ge 1$ sufficiently large} since this contradicts \eqref{basicregpropviol} and \eqref{basicregpropvio2}.
\end{itemize}

\subsection{Indices of discontinuity}\label{discontsect}
The following explanations illustrate the difference between $\phi=+\infty$, yielding the (log-)convex minorant, and $\phi\neq +\infty$ where so-called \emph{indices of discontinuity} can appear. So let now $\phi\neq +\infty$:

\begin{itemize}
\item[$(a)$] Take a principal index $p_i$, let $t_1<t_2$ and assume that $S_{p_i}\in D_{t_1}$ and that $S_{p_i}\in D_{t_2}$. Then $S_{p_i}\in D_t$ holds for all $t\in[t_1,t_2]$ as well. So for each $p_i$ we can associate an interval $I_i\subseteq\RR$ such that $S_{p_i}\in D_t$ for all $t\in I_i$.

\item[$(b)$] We gather some properties for these intervals:

\begin{itemize}
\item[$(b.1)$] For each $i\ge 1$ the length of $I_i$ is finite because otherwise this contradicts $a_p=+\infty$ for possibly at most only finitely many $p\ge 1$.

\item[$(b.2)$] $I_i$ is, in general, of the form $I_i=[t_i,\tau_i)$ and $t_i\in I_i$ holds because $\phi$ is continuous. For $i=0$ we have $I_0=[-\infty,\tau_0)$ and the slope $-\infty$ is formally corresponding to the $y$-axis on which only $S_0$ is lying.

\item[$(b.3)$] Clearly $\tau_i=t_{i+1}$ is valid for all $i\in\NN$.

\item[$(b.4)$] One might have $I_i=t_i=t_{i+1}$: If $t_i=t_{i+1}=\dots=t_{i+k}$ for some $k\in\NN_{>0}$, then all points $S_{p_i},\dots,S_{p_{i+k}}$ are lying on the joint line segment $D_{t_i}$ with slope $t_i$. However, this cannot be the case for $I_0$ since then $S_{p_1}$ would be lying on the $y$-axis.

\item[$(b.5)$] Finally, $t_i\rightarrow+\infty$ as $i\rightarrow+\infty$ is also clear: $\phi(t)\rightarrow+\infty$ as $t\rightarrow+\infty$ but for each $t\in\RR$ we only deal with finitely many $S_p\in B_t$ and since $a_p=+\infty$ for only possibly finitely many $p\ge 1$.
\end{itemize}

\item[$(c)$] Let $p_i$ be a principal index. If $\tau_i\notin I_i$ or if $\tau_i\in I_i$ depends on the following relation between $\phi(\tau_i)$ and $p_{i+1}$:
    \begin{itemize}
    \item[$(c.1)$] The case $\phi(\tau_i)=\phi(t_{i+1})<p_{i+1}$ is not possible by continuity of $\phi$, since $\phi$ is non-decreasing and by definition of $I_i$ and $p_{i+1}$: In this situation the requirement $S_{p_{i+1}}\in D_{\tau_i}=D_{t_{i+1}}$ is impossible.

    \item[$(c.2)$] If $\phi(\tau_i)=\phi(t_{i+1})>p_{i+1}$, then $I_i=[t_i,t_{i+1}]$ since $S_{p_{i+1}}\in B_{t_{i+1}}$.

    \item[$(c.3)$] If $\phi(\tau_i)=\phi(t_{i+1})=p_{i+1}$, then $I_i=[t_i,t_{i+1}]$ if and only if the points $S_{p_i}$, $S_{p_{i+1}}$ both lie on the segment $D_{t_{i+1}}$. However, this fact is \emph{not clear in general.}

         Principal indices $p_{i+1}$ with this property are called \emph{indices of discontinuity.}

    \item[$(c.4)$] We illustrate this situation: Assume that $a_{p_i}=a_{p_{{i+1}}}$, and so the points $S_{p_i}$, $S_{p_{i+1}}$ have the same $y$-coordinate, and that all $S_p$, $0\le p<p_{i+1}$ with $p\neq p_i$, have the property that they are lying strictly above each segment $D_t$, $0<t_i\le t<t_{i+1}$, which is passing through $S_{p_i}$ and having positive slope $t$.

        Then, by definition, for $t=t_{i+1}$ the segment $D_{t}$ is passing through $S_{p_{i+1}}$ but $S_{p_i}$ is lying strictly above $D_{t_{i+1}}$ and so $p_{i+1}$ is an index of discontinuity.
    \end{itemize}

\item[$(d)$] Recall: $t_i=t_{i+1}=\dots=t_{i+k}$ for some $k\in\NN_{>0}$ means that all points $S_{p_i},\dots,S_{p_{i+k}}$ are lying on a joint line segment and hence $p_{i+1},\dots,p_{i+k}$ are not indices of discontinuity.

    Note also that $p_0=0$ is never an index of discontinuity.
\end{itemize}

\vspace{6pt}
Now let us consider the case where $\phi(T)=+\infty$ for some $T\in\RR$. At the beginning, for all sufficiently small $t$ we can proceed as above for $\phi\neq +\infty$. Note that there exist at least finitely many principal indices $p_i$; indeed there exists at least one such index, namely $p_0=0$. As $t\rightarrow T$ we have that $\phi(t)\rightarrow+\infty$, i.e. the width of $B_t$ is tending to $+\infty$, but $D_t\rightarrow D_T$ and so the slopes of the segments $D_t$ under consideration are uniformly bounded by $T$. Thus we have to distinguish; see $(iv)$ in Section \ref{principalsect}:

\begin{itemize}
\item[$(a)$] There exist finitely many principal indices: First, since the set of indices of discontinuity is a subset of all principal indices there exist for sure also at most finitely many indices of discontinuity. When $p_{i_0}$ is denoting the last principal index, then we set $I_{i_0}:=[t_{i_0},T)$.

\item[$(b)$] There exist infinitely many principal indices: Here the situation is analogous to the study of $\phi\neq +\infty$ before. However, note that the length of the intervals $I_i$ is tending to $0$ as $i\rightarrow+\infty$ since the maximal admissible slope is bounded by $T$.
\end{itemize}

Finally, let us consider $\phi=+\infty$:

\begin{itemize}
\item[$(a)$] In any case, the set of indices of discontinuity is empty since $(c.2)$ before holds automatically for all $i\in\NN$; see Sections \ref{lcminorsect} and \ref{nonstandardsect}.

\item[$(b)$] More precisely, in the exotic situation treated in Section \ref{case1sect} we have $I_0=\{-\infty\}$ and this is the only ``interval'' under consideration.

\item[$(c)$] In the case treated in Section \ref{case2sect} we have to distinguish: If there exist infinitely many principal indices, then to each of these we can assign an interval $I_i=[t_i,t_{i+1}]$ and with $I_0=[-\infty,t_1]$. Here the length of $I_i$ is tending to $0$ since the maximal admissible slope is bounded by $\mathbf{a}_{\iota}$. If there exist only finitely many principal indices, say $i_0$ many, then the last interval is given by $I_{i_0}=[t_{i_0},\mathbf{a}_{\iota})$.

    \emph{Note:} In the case of $(i)$ in Example \ref{nonstandardsectex} we have $p_0=0$ and $p_1=1$ and this is the last principal index. Then $I_0=[-\infty,-1]$ and $I_1=[-1,c)$ and no further intervals appear. In $(ii)$ there we have $I_0=[-\infty,c)$ and no further interval is considered.
\end{itemize}

\subsection{Definition of the regularized sequence}\label{defregsect}
We define the crucial sequence of line segments $(L_i)_{i\in\NN}$ as follows; first we focus on $\phi\neq +\infty$ and comment then on the other situations:

\begin{itemize}
\item[$(I)$] If $p_i$ and $p_{i+1}$ are principal indices, then $L_i$ shall denote the segment of the straight line on $[p_i,p_{i+1})$, i.e. $\pr_1(L_i)=[p_i,p_{i+1})$, such that $L_i$ is having slope $\tau_i=t_{i+1}$ and it is passing through the point $S_{p_i}=(p_i,a_{p_i})$. Thus, we get
    $$L_i: s\mapsto t_{i+1}s+a_{p_i}-t_{i+1}p_i=t_{i+1}(s-p_i)+a_{p_i},\;\;\;p_i\le s<p_{i+1},\;i\in\NN.$$

\item[$(II)$] If $p_{i+1}$ is not an index of discontinuity, then we have $L_i\subseteq D_{t_{i+1}}$. If $p_{i+1}$ is an index of discontinuity, then $L_i$ is parallel to $D_{t_{i+1}}$ but is strictly lying over $D_{t_{i+1}}$ because $\phi(t_{i+1})=\phi(\tau_i)=p_{i+1}$ and $S_{p_i}$ lies strictly above $D_{t_{i+1}}$. The slopes of $L_i$ are non-decreasing and tending to $+\infty$ as $i\rightarrow+\infty$.
\end{itemize}

Then consider
\begin{equation}\label{setB}
\mathcal{L}:=\bigcup_{i\in\NN}L_i,
\end{equation}
and hence $\mathcal{L}\subseteq\RR^2$ consists of a set of straight lines having the following properties:

\begin{itemize}
\item[$(*)$] At all principal indices, except at $p_0=0$, cusps appear and at all indices of discontinuity a jump.

\item[$(*)$] More precisely, between $p_0$ and the first index of discontinuity, and then between two consecutive indices of discontinuity, the set of lines is connected, having cusps at all principal indices. Finally it is closed on the left-hand side and open on the right-hand side at the jumps.
\end{itemize}

We introduce the regularizing mapping $r^{\phi}$ by
\begin{equation}\label{regrelation}
S_p=(p,a_p)\mapsto S_p^{\phi}:=(p,a_p^{\phi}),
\end{equation}
where $S_p^{\phi}$ denotes the \emph{vertical projection} of the point $S_p$ onto the set $\mathcal{L}$; i.e. $a_p$ is projected vertically down onto $a_p^{\phi}$. This gives the regularized set $\mathcal{S}^{\phi}:=\{S_p^{\phi}=(p,a_p^{\phi}): p\in\NN\}$ and the regularized sequence $\mathbf{a}^{\phi}:=(a^{\phi}_p)_{p\in\NN}$ satisfies
\begin{itemize}
\item[$(\alpha)$] $a^{\phi}_p\le a_p$ for all $p\in\NN$,

\item[$(\beta)$] $(\mathbf{a}^{\phi})^{\phi}=\mathbf{a}^{\phi}$,

\item[$(\gamma)$] $a^{\phi}_p=a_p$ if and only if $p$ is a principal index. In particular, we have $a^{\phi}_0=a_0$.
\end{itemize}

We comment now on the other cases for $\phi$. The corresponding sequence $S_p^{\phi}:=(p,a_p^{\phi})$ is defined analogously via \eqref{regrelation} but the set $\mathcal{L}$ can differ from the case $\phi\neq +\infty$.

\begin{itemize}
\item[$(a)$] Let $\phi(T)=+\infty$ for some $T\in\RR$ and distinguish:
\begin{itemize}
\item[$(a.1)$] If there exist infinitely many principal indices, then the only difference compared with $\phi\neq +\infty$ is that the slopes of the lines $L_i$ are tending to $T$ as $i\rightarrow+\infty$.

\item[$(a.2)$] If there exist finitely many principal indices, say $i_0$ many, then the set $\mathcal{L}$ in \eqref{setB} has the form
$$\mathcal{L}:=\bigcup_{0\le i\le i_0}L_i;$$
i.e. it consists only of $i_0+1$ lines. The first $i_0$ lines are as described before and the ``last member'' $L_{i_0}$ is defined to be the line with slope $T$ and passing through the point $S_{p_{i_0}}$; i.e. we get $L_{i_0}:=D_{T}^{-}:=\lim_{\varepsilon\rightarrow 0}D_{T-\varepsilon}$.

In view of this definition we see that in this case formally it makes sense to consider $\mathbf{a}^{\phi}$ even if $\mathbf{a}=(a_p)_{p\in\NN}$ satisfies $a_p=+\infty$ for all $p\ge 1$ sufficiently large: One has only finitely many principal indices and \eqref{regrelation} makes sense since for all $p\in\NN$ large the $y$-coordinate $+\infty$ is ``projected vertically down'' onto the line $L_{i_0}$.
\end{itemize}

\item[$(b)$] Let $\phi=+\infty$ and distinguish:

\begin{itemize}
\item[$(b.1)$] If $\mathbf{a}$ satisfies \eqref{basicregprop}, then $\mathcal{L}$ is precisely the set of lines described in Section \ref{lcminorsect}. More precisely, $\mathcal{L}=\bigcup_{i\in\NN}L_i$ is connected since the set of indices of discontinuity is empty and all further properties of $\mathcal{L}$ are the same as for $\phi<+\infty$.

\item[$(b.2)$] In the exotic situation treated in Section \ref{case1sect} formally we have $\mathcal{L}=\{L_0\}$ and $L_0$ denotes the $y$-axis since this is the only line under consideration. The set of principal indices is equal to $\{p_0\}$ and the set of indices of discontinuity is empty. However, in order to make the definition \eqref{regrelation} fit with the approach in Section \ref{case1sect} we put $\mathcal{L}=\{L_0,L_1\}$ with $L_1: y=-\infty$ since all $S_p$, $p\ge 1$, are formally projected down onto this line.

\item[$(b.3)$] In the situation of Section \ref{case2sect} the set of indices of discontinuity is again empty. However, we have either $\mathcal{L}=\bigcup_{i\in\NN}L_i$, if there exist infinitely many principal indices, or $\mathcal{L}:=\bigcup_{0\le i\le i_0}L_i$ else. In the latter case the last member $L_{i_0}$ is denoting the line passing through the ``last principal point'' $S_{p_{i_0}}$ and having slope $\mathbf{a}_{\iota}$; i.e. $L_{i_0}$ coincides with $D_{\mathbf{a}_{\iota}}^{-}:=\lim_{\varepsilon\rightarrow 0}D_{\mathbf{a}_{\iota}-\varepsilon}$ on $[p_{i_0},+\infty)$.
\end{itemize}
\end{itemize}

In all instances properties $(\alpha)$, $(\beta)$, $(\gamma)$ hold. However, for $(\beta$) in $(b.2)$ recall the convention there, so note that $\mathbf{a}^{\phi}\notin\RR^{\NN}$ and thus $\mathbf{a}\mapsto\mathbf{a}^{\phi}$ is not a mapping $\RR^{\NN}\rightarrow\RR^{\NN}$ anymore.

\subsection{Comparison between regularizing functions}
Let $\mathbf{a}\in\RR^{\NN}$ be fixed.

\begin{itemize}
\item[$(*)$] Given two regularizing functions $\phi_1$ and $\phi_2$ with $\phi_1\le\phi_2$ then we obtain immediately $\mathbf{a}^{\phi_2}\le\mathbf{a}^{\phi_1}$.

\item[$(*)$] Since the ``largest possible'' or ``most extreme'' example is given by $\phi=+\infty$ we have
    \begin{equation}\label{comparisonregequ}
    \mathbf{a}^c\le\mathbf{a}^{\phi}\le\mathbf{a},
    \end{equation}
    for all $\mathbf{a}\in\RR^{\NN}$ and all regularizing functions $\phi$.

\item[$(*)$] In view of \eqref{comparisonregequ} the constructions treated in Sections \ref{lcminorsect} and \ref{nonstandardsect} yield the smallest possible sequence within this geometric approach. Indeed, in Section \ref{case1sect} one obtains the smallest possible regularized ``sequence''.
\end{itemize}

\subsection{The counting function}\label{countfctsect}
Let $\mathbf{a}\in\RR^{\NN}$ and $\phi$ be a regularizing function. Then set
\begin{equation}\label{functionm}
m^{\phi}(t):=\sup_{p\in\NN}\{S_p\in D_t\},\;\;\;t\in\{-\infty\}\cup J^{\phi},
\end{equation}
where $J^{\phi}\subseteq\RR$ is the natural ``interval of definition'' depending on $\phi$. Recall that $D_{-\infty}$ corresponds to the $y$-axis and let $J^{\phi}_r$ be the right-end point of $J^{\phi}$; i.e. the limit point. We describe this situation in detail:

\begin{itemize}
\item[$(a)$] For any regularizing function $\phi$ one has $m^{\phi}(t)=0$ for all sufficiently small $t$ since in this case $S_0$ is the only point lying on $D_t$. We also put $m^{\phi}(-\infty):=0$ since $S_0$ is clearly the only point belonging to $\mathcal{S}$ and lying on the $y$-axis.

\item[$(b)$] If $\phi\neq +\infty$, then $m^{\phi}:\{-\infty\}\cup\RR\rightarrow\NN_{>0}$ is well-defined; i.e. $m^{\phi}(t)<+\infty$ for all $t\in\{-\infty\}\cup\RR$, because there are only finitely many points $S_p$ which belong to $B_t$ and hence are possibly lying on the segment $D_t$. Thus we have $J^{\phi}=\RR$ and $J^{\phi}_r=+\infty$.

\item[$(c)$] Similarly, if $\phi(T)=+\infty$ for some $T\in\RR$, then $m^{\phi}$ is well-defined on $[-\infty,T)$ since for each $t\in[-\infty,T)$ only finitely many points $S_p$ can belong to $B_t$. In this case $J^{\phi}=(-\infty,T)$ and so $J^{\phi}_r=T$.

\item[$(d)$] Let $\phi=+\infty$, then we distinguish:

\begin{itemize}
\item[$(d.1)$] $m^{\phi}$ is well-defined on $\{-\infty\}\cup\RR$ provided that $\mathbf{a}$ satisfies \eqref{basicregprop}: When $t\in\RR$ is fixed, then eventually all $S_p$ are lying strictly above the half-ray $D_t$; see Section \ref{normgeomconstrsect}. Hence $J^{\phi}=\RR$ and $J^{\phi}_r=+\infty$ is valid.

\item[$(d.2)$] If \eqref{basicregprop} is violated and $\mathbf{a}_{\iota}>-\infty$, see \eqref{basicregpropvio2} and \eqref{aiotadef}, then $m^{\phi}$ is well-defined on $[-\infty,\mathbf{a}_{\iota})$ and so $J^{\phi}=(-\infty,\mathbf{a}_{\iota})$ and $J^{\phi}_r=\mathbf{a}_{\iota}$; see Section \ref{case2sect}. However, note that the value $\mathbf{a}_{\iota}$ is excluded since in general one might have $m^{\phi}(\mathbf{a}_{\iota})=+\infty$; see $(ii)$ in Example \ref{nonstandardsectex}:

    If $c\in\RR$ is fixed and $a_p:=cp$ for all $p\in\NN$, then $m^{\phi}(t)=0$ for $t\in[-\infty,c)$ and $m^{\phi}(c)=+\infty$. Recall that this sequence is already convex and all $S_p$ are lying on the straight line passing through the origin and having slope $c$. This is a situation which is excluded by the standard requirement \eqref{basicregprop}.

\item[$(d.3)$] If $\mathbf{a}_{\iota}=-\infty$, i.e. \eqref{basicregpropviol} in Section \ref{case1sect}, then we consider $m^{\phi}$ only at $-\infty$ with $m^{\phi}(-\infty):=0$. In this case we have $J^{\phi}=\emptyset$.
\end{itemize}
\end{itemize}

Note that in all cases clearly $\lim_{t\rightarrow J^{\phi}_r}\phi(t)=+\infty$ and for $(d.3)$ we have $\phi=+\infty$ in any case. Next we summarize some properties for $m^{\phi}$:

\begin{itemize}
\item[$(x)$] $m^{\phi}$ is non-decreasing on $J^{\phi}$; for case $(d.3)$ this statement is empty: Let $t_1<t_2$, then we distinguish two cases. If $S_{m^{\phi}(t_1)}\in D_{t_2}$, then obviously $m^{\phi}(t_1)\le m^{\phi}(t_2)$ holds. If $S_{m^{\phi}(t_1)}$ lies above the segment $D_{t_2}$, then all points $S_p$ with $0\le p\le m^{\phi}(t_1)$ lie above $D_{t_2}$ as well. Thus we obtain in this case even $m^{\phi}(t_1)<m^{\phi}(t_2)$.

\item[$(xx)$] Now we study the limit case; i.e. the behavior of $m^{\phi}$ at the right end point of $J^{\phi}$ again except the trivial ``constant/one-point'' case in $(d.3)$.

\begin{itemize}
\item[$(xx.1)$] If $\phi\neq +\infty$, then $\lim_{t\rightarrow+\infty}m^{\phi}(t)=+\infty$: Suppose that there exists $k\in\NN$ such that $m^{\phi}(t)\le k$ for all $t\in\{-\infty\}\cup\RR$, then $a_p=+\infty$ for all $p>k$ follows which is a contradiction to the basic assumption on $\mathbf{a}$.

\item[$(xx.2)$] If $\phi=+\infty$ and such that $\mathbf{a}$ satisfies \eqref{basicregprop}, then $\lim_{t\rightarrow+\infty}m^{\phi}(t)=+\infty$ holds by similar reasons; see Section \ref{normgeomconstrsect}.

\item[$(xx.3)$] If $\phi=+\infty$, \eqref{basicregprop} is violated and $\mathbf{a}_{\iota}>-\infty$, see case $(d.2)$ before, then we have that $\lim_{t\rightarrow\mathbf{a}_{\iota}}m^{\phi}(t)=+\infty$ if and only if there exist infinitely many principal indices; see Example \ref{nonstandardsectex1} and $(iii)$ in Section \ref{principalsect}.

    So here this limit might be finite and for this consider Example \ref{nonstandardsectex}: The sequence given by \eqref{nonstandardsectexequ} in $(i)$ yields $\lim_{t\rightarrow\mathbf{a}_{\iota}}m^{\phi}(t)=1$ and in $(ii)$ we even get $\lim_{t\rightarrow\mathbf{a}_{\iota}}m^{\phi}(t)=0$.

\item[$(xx.4)$] If $\phi(T)=+\infty$ for some $T\in\RR$, then similarly the value of $\lim_{t\rightarrow T}m^{\phi}(t)$ depends on the existence of finitely or infinitely many principal indices and so on the relation between $T$ and $\mathbf{a}$. If $\mathbf{a}$ satisfies \eqref{liminfT}, then $\lim_{t\rightarrow T}m^{\phi}(t)=+\infty$, but if $\mathbf{a}$ eventually lies over lines having slope $T$, then this value is finite and can even be equal to $0$.
\end{itemize}

\item[$(xxx)$] $m^{\phi}$ is continuous from the right on $J^{\phi}$. Indeed, in the case treated in $(d.3)$ there is nothing to show. For the segments $D_t$ we have $D_t=D_t^{+}$, where
$$D_t^{+}:=\lim_{\varepsilon\rightarrow 0}D_{t+\varepsilon},$$
hence the mapping $t\mapsto D_t$ is right-continuous. If $p_i$ is not an index of discontinuity, then we obtain $D_t=D_t^{-}$ and recall the notation $D_t^{-}:=\lim_{\varepsilon\rightarrow 0}D_{t-\varepsilon}$. Analogously we define $(m^{\phi})^{+}(t)$ and with this notation we have to show $m^{\phi}(t)=(m^{\phi})^{+}(t)$ for all $t\in J^{\phi}$. Because $m^{\phi}(t)\in\NN$ for all $t\in J^{\phi}$ we have $(m^{\phi})^{+}(t)\in\NN$. Thus there exists $\varepsilon_0>0$ such that $m^{\phi}(t+\varepsilon)=(m^{\phi})^{+}(t)$ for all $\varepsilon>0$ with $\varepsilon\le\varepsilon_0$. Then $S_{(m^{\phi})^{+}(t)}\in D_{t+\varepsilon}$ holds for all those $\varepsilon$, hence $S_{(m^{\phi})^{+}(t)}\in D_t$ because $D_t=D_t^{+}$. Finally $(m^{\phi})^{+}(t)\le m^{\phi}(t)$ implies the desired property.
\end{itemize}

\subsection{The trace function}\label{tracefctsct}
For each $\mathbf{a}\in\RR^{\NN}$ and each regularizing function $\phi$ we can associate a weight function, the so-called \emph{trace function} $A^{\phi}: J^{\phi}\rightarrow\RR$. In order to lighten notation we suppress the dependence on the given sequence $\mathbf{a}$ except in \eqref{tracefctstable}. $A^{\phi}$ is defined analogously as in \eqref{tracefunction} as follows:\vspace{6pt}

For $t\in J^{\phi}$ the value $-A^{\phi}(t)$ denotes the $y$-coordinate of the unique intersection point of the segment $D_t$ with the line $\{(0,y): y\in\RR\}$; i.e. with the $y$-axis. Moreover, we extend the definition and put $A(-\infty):=-a_0$ and so $A^{\phi}: \{-\infty\}\cup J^{\phi}\rightarrow\RR$.

If we take an arbitrary straight line $\ell$ with slope $t\in\RR$ and which is passing through some point $S_p=(p,a_p)$ with $S_p\in B_t$, then the $y$-coordinate of the intersection point of $\ell$ and the $y$-axis is equal to $a_p-pt$.

Let us gather some immediate consequences:

\begin{itemize}
\item[$(*)$] If $J^{\phi}=\emptyset$, i.e. in the situation of case $(d.3)$ in Section \ref{countfctsect}, then $A(-\infty):=-a_0$ and $-\infty$ the only value for which $A$ is defined.
\end{itemize}

Assume now that $J^{\phi}$ is non-empty.

\begin{itemize}
\item[$(*)$] By definition $A^{\phi}$ is strictly increasing on $J^{\phi}$.

\item[$(*)$] $\lim_{t\rightarrow-\infty}A^{\phi}(t)=-a_0$ holds since for all $t$ sufficiently small we have that $S_0=(0,a_0)$ is the only point belonging to $B_t$; for this recall that $\lim_{t\rightarrow-\infty}\phi(t)=0$ and so $m^{\phi}(t)=0$ for all $t$ sufficiently small. The same is valid if $\phi=+\infty$ for the cases $(d.1)$ and $(d.2)$ in Section \ref{countfctsect},

\item[$(*)$] $A^{\phi}$ is right-continuous and it is not continuous at $t=t_i$ if $p_i$ is an index of discontinuity yielding a jump. Consequently, for $\phi=+\infty$ the function $A^{\phi}$ is continuous.

\item[$(*)$] $\lim_{t\rightarrow+\infty}A^{\phi}(t)=+\infty$ holds if the right-end point of $J^{\phi}$ is equal to $+\infty$, so if $\phi<+\infty$ or if $\phi=+\infty$ and $\mathbf{a}$ has \eqref{basicregprop}.

More precisely, in these cases we are allowed to consider arbitrary large slopes $t$ and therefore the $y$-coordinates $a_p-pt$ cannot be uniformly bounded away from $-\infty$ for all $p\in\NN_{>0}$ as $t\rightarrow+\infty$. Otherwise, for all $t$ large enough at least finitely many points $S_p$, with $1\le p\le p_0$ for some $p_0\in\NN_{>0}$ depending on this fixed $t$, will lie strictly below $D_t$.

\item[$(*)$] If $\phi(T)=+\infty$ for some $T\in\RR$, then $\lim_{t\rightarrow J^{\phi}_r}A^{\phi}(t)=\lim_{t\rightarrow T}A^{\phi}(t)$ can be finite or not depending on the growth behavior of $\mathbf{a}$. For example one has $\lim_{t\rightarrow T}A^{\phi}(t)<+\infty$ provided there exist finitely many principal indices: Let $p_{i_0}$ be the last one, then $\lim_{t\rightarrow T}A^{\phi}(t)=-(a_{p_{i_0}}-Tp_{i_0})=Tp_{i_0}-a_{p_{i_0}}$; see \eqref{lineT} in $(iv)$ in Section \ref{principalsect}.

\item[$(*)$] Moreover, in the situation of $(d.2)$ in Section \ref{countfctsect}, i.e. if $\phi=+\infty$, \eqref{basicregprop} is violated and $\mathbf{a}_{\iota}>-\infty$, then we also have that $\lim_{t\rightarrow J^{\phi}_r}A^{\phi}(t)=\lim_{t\rightarrow\mathbf{a}_{\iota}}A^{\phi}(t)<+\infty$: The admissible slope of the lines under consideration is bounded by $\mathbf{a}_{\iota}$ and the limiting case is a straight line with slope $\mathbf{a}_{\iota}$ and being parallel to the line $t\mapsto t\mathbf{a}_{\iota}$; see $(i)$ in Example \ref{nonstandardsectex}. The negative value of the $y$-coordinate of the intersection point of this line with the $y$-axis is the limit under consideration. If $a_p=cp$ for $c\in\RR$ and all $p\in\NN$, then in fact $A(t)=0$ for all $t\in(-\infty,c)$ and even in the limit we get $A(c)=0$; see $(ii)$ in Example \ref{nonstandardsectex}.
\end{itemize}

We comment in more detail on the particular cases:

\begin{itemize}
\item[$(i)$] If $\phi\neq +\infty$, then by definition of $D_t$ and the trace function we see that $a_p-pt\le-A^{\phi}(t)\Leftrightarrow A^{\phi}(t)\le p t-a_p$ for all $0\le p\le\phi(t)$ and equality holds, in particular, for $p=m^{\phi}(t)$. Thus we summarize, compare with \eqref{trace}:

\begin{equation}\label{Regulariz2}
A^{\phi}(t)=\max_{p\in\NN: \;0\le p\le\phi(t)}\{pt-a_p\}=m^{\phi}(t)t-a_{m^{\phi}(t)},\;\;\;t\in J^{\phi}=\RR.
\end{equation}

\item[$(ii)$] If $\phi(T)=+\infty$ for some $T\in\RR$, then \eqref{Regulariz2} holds for all $t\in(-\infty,T)$.

\item[$(iii)$] Let $\phi=+\infty$, then we distinguish:

\begin{itemize}
\item[$(iii.1)$] If \eqref{basicregprop} holds, then \eqref{Regulariz2} is valid for all $t\in\RR(=J^{\phi})$ and $\max_{p\in\NN: \;0\le p\le\phi(t)}$ is replaced by $\sup_{p\in\NN}$.

\item[$(iii.2)$] If \eqref{basicregprop} is violated and $\mathbf{a}_{\iota}>-\infty$, then \eqref{Regulariz2} follows analogously for all $t\in J^{\phi}=(-\infty,\mathbf{a}_{\iota})$ and $\max_{p\in\NN: \;0\le p\le\phi(t)}$ is replaced again by $\sup_{p\in\NN}$.

\item[$(iii.3)$] Finally, if \eqref{basicregprop} is violated and $\mathbf{a}_{\iota}=-\infty$, i.e. \eqref{basicregpropviol} in Section \ref{case1sect}), then we have $J^{\phi}=\emptyset$ and recall that we put $A^{\phi}(-\infty):=-a_0=-a_{m^{\phi}(-\infty)}$ and $-\infty$ is the only point at which formally $A^{\phi}$ is defined.
\end{itemize}
\end{itemize}
Summarizing, $m^{\phi}$ and $A^{\phi}$ are related by
\begin{equation}\label{Regulariz2new}
A^{\phi}(t)=tm^{\phi}(t)-a_{m^{\phi}(t)},\;\;\;t\in\;J^{\phi};
\end{equation}
i.e. $m^{\phi}$ is the right-derivative of $A^{\phi}$. Indeed, with $A(-\infty):=-a_0$, see the comment in $(iii.3)$, and the convention $0\cdot(-\infty):=0$ we have that \eqref{Regulariz2new} can be extended to $t=-\infty$ in all instances.\vspace{6pt}

To each sequence $\mathbf{a}\in\RR^{\NN}$ and each regularizing function $\phi$ we have constructed a unique trace function $A^{\phi}$. But there exist different sequences which lead to the same trace function. In particular, and this is immediate by the definitions given in Section \ref{defregsect}, we have that $\mathbf{a}$ and $\mathbf{a}^{\phi}$ yield the same trace function, i.e.
\begin{equation}\label{tracefctstable}
A_{\mathbf{a}^{\phi}}^{\phi}=A_{\mathbf{a}}^{\phi},
\end{equation}
for all $\mathbf{a}\in\RR^{\NN}$ and all regularizing functions $\phi$. This fact even holds in the exotic situation $(iii.3)$ since $A^{\phi}_{\mathbf{a}}(-\infty)=-a_0=-a^{\phi}_0=A_{\mathbf{a}^{\phi}}^{\phi}(-\infty)$ and $-\infty$ is the only point for which \eqref{tracefctstable} should be verified but, formally, \eqref{tracefctstable} can be extended to whole $\RR$ by setting $A_{\mathbf{a}^{\phi}}^{\phi}(t)=+\infty=A_{\mathbf{a}}(t)$ for all $t\in\RR$.

\subsection{Computing the regularized sequence via the trace function}\label{regsequfromtracesect}
In this final Section we show how, when $\mathbf{a}$ is given, the trace function $A^{\phi}\equiv A_{\mathbf{a}}^{\phi}$ can be used to compute $\mathbf{a}^{\phi}$ and hence $\mathcal{S}^{\phi}$. Summarizing, one is able to study $\mathbf{a}\mapsto A_{\mathbf{a}}^{\phi}\mapsto\mathbf{a}^{\phi}$.

\begin{theorem}\label{Regulariz3}
Let $\mathbf{a}\in\RR^{\NN}$ and a regularizing function $\phi$ be given with $J^{\phi}\neq\emptyset$. Then, the sequence $\mathbf{a}^{\phi}=(a^{\phi}_p)_{p\in\NN}$ can be expressed by
\begin{equation}\label{Regulariz3-equ1}
a^{\phi}_p=\sup_{t\in J^{\phi}: \phi(t)\ge p}\{pt-A^{\phi}(t)\},\;\;\;p\in\NN.
\end{equation}
Therefore recall that $\lim_{t\rightarrow J^{\phi}_r}\phi(t)=+\infty$ and so each $p\in\NN$ can be computed via \eqref{Regulariz3-equ1}.
\end{theorem}

For the particular cases we mention:

\begin{itemize}
\item[$(i)$] If $J^{\phi}=\emptyset$, i.e. if $\phi=+\infty$ and $\mathbf{a}$ has \eqref{basicregpropviol}, then \eqref{Regulariz3-equ1} is formally not applicable; see also Section \ref{case1sect}. But when using the conventions $0\cdot(-\infty):=0$ and $p\cdot(-\infty):=-\infty$ for any $p\in\NN_{>0}$ and recalling $A^{\phi}(-\infty):=-a_0$, then in \eqref{Regulariz3-equ1} one can consider in any case $t\in\{-\infty\}\cup J^{\phi}$. And this gives here again $a^{\phi}_0=a_0$ and $a^{\phi}_p=-\infty$ for $p\ge 1$.

\item[$(ii)$] If $\phi=+\infty$ and $\mathbf{a}$ has \eqref{basicregprop}, then $\sup_{t\in J^{\phi}: \phi(t)\ge p}$ can be replaced by $\sup_{t\in\RR}$.

\item[$(iii)$] If $\phi=+\infty$ and $\mathbf{a}$ has \eqref{basicregpropvio2}, then $\sup_{t\in J^{\phi}: \phi(t)\ge p}$ can be replaced by $\sup_{t\in(-\infty,\mathbf{a}_{\iota})}$.

\item[$(iv)$] If $\phi\neq+\infty$ resp. $\phi(T)=+\infty$ for some $T\in\RR$, then $\sup_{t\in J^{\phi}: \phi(t)\ge p}$ can be replaced by $\sup_{t\in\RR: \phi(t)\ge p}$ resp. $\sup_{t\in(-\infty,T): \phi(t)\ge p}$. However, in general we cannot consider all $t\in\RR$ resp. $t\in(-\infty,T)$ since in the regularization procedure segments $D_t$ with $t$ such that $\phi(t)<p$ are not considered; see $(c.4)$ in Section \ref{discontsect} and the arguments in $(a)$ resp. $(c)$ in the proof below.

\item[$(v)$] When $J^{\phi}\neq\RR$, then let us set $(J^{\phi})^c:=\RR\backslash J^{\phi}$ and, analogously as the conventions used in Sections \ref{case1sect} and \ref{case2sect}, we put $A^{\phi}(t):=+\infty$ for all $t\in(J^{\phi})^c$. This allows then to replace $\sup_{t\in(-\infty,\mathbf{a}_{\iota})}$ resp. $\sup_{t\in(-\infty,T): \phi(t)\ge p}$ by $\sup_{t\in\RR}$ resp. $\sup_{t\in\RR: \phi(t)\ge p}$; for the second situation recall that we put $\phi(t)=+\infty$ for all $t\ge T$. Finally, in case $(i)$ we can also put $A^{\phi}(k):=+\infty$ for all $k\in\RR$ and so consider $\sup_{\{-\infty\}\cup\RR}$ instead; see again Section \ref{case1sect}.
\end{itemize}

\demo{Proof}
We distinguish between the different cases.

\begin{itemize}
\item[$(a)$] Assume that $\phi\neq +\infty$. Let $p\in\NN$ and consider the auxiliary line $\ell_p:=\{(x,y)\in\RR^2: x=p\}$. For all $t\in J^{\phi}=\RR$ with $\phi(t)\ge p$ consider the segment $D_t$. For $p=0$ the line $\ell_p$ is the $y$-axis and we have $\phi(t)\ge 0$ for any $t\in\RR$. Hence $D_t$ intersects $\ell_p$ in $(p,pt-A^{\phi}(t))$ but, by definition, the point $S_p^{\phi}=(p,a_p^{\phi})$ is not lying below $D_t$ which implies $pt-A^{\phi}(t)\le a_p^{\phi}$. Note: For $p=0$ we get $-A^{\phi}(t)\le a_0=a_p^{\phi}$ for all $t\in\RR$ and equality is attained in the limit $t=-\infty$.

    For the converse, first by construction resp. definition of $\mathbf{a}^{\phi}$ given in Section \ref{defregsect} we see that for all $p\in\NN$ there exists $i\in\NN$ such that $S_p^{\phi}\in L_i$ and that the desired equality holds for any principal index since $a_{p_i}=a_{p_i}^{\phi}$. Let $p\in\NN$ and $p_{i+1}$ be now the smallest principal index satisfying $p_{i+1}>p>p_i$. Indeed, recall that by $(ii)$ in Section \ref{principalsect} there exist infinitely many principal indices.

    We distinguish now between two cases, see $(c)$ in Section \ref{discontsect}: If the index $p_{i+1}$ is not an index of discontinuity, then we have $L_i\subseteq D_{t_{i+1}}$ and because $\phi(t_{i+1})=\phi(\tau_i)\ge p_{i+1}>p$ holds we are done since the supremum in \eqref{Regulariz3-equ1} is attained at $t=t_{i+1}$. If $p_{i+1}$ is an index of discontinuity, and so necessarily $\phi(\tau_i)=\phi(t_{i+1})=p_{i+1}$, then $L_i\subseteq D_{t_{i+1}}^{-}$ where $D_{t_{i+1}}^{-}:=\lim_{\varepsilon\rightarrow 0}D_{t_{i+1}-\varepsilon}$. Thus \eqref{Regulariz3-equ1} holds as well but the supremum is not attained for $t\in\RR=J^{\phi}$.

\item[$(b)$] Assume that $\phi=+\infty$ and distinguish:
\begin{itemize}
\item[$(b.1)$] If $\mathbf{a}$ satisfies \eqref{basicregprop}, then the same arguments from $(a)$ apply and note that the set of indices of discontinuity is empty.

\item[$(b.2)$] Let $\mathbf{a}$ have \eqref{basicregpropvio2}. Again the set of indices of discontinuity is empty. If there exist infinitely many principal indices, then the result follows by repeating the arguments from $(a)$. If there exist only finitely many principal indices, say $i_0$ many, then for all $p<p_{i_0}$ we can proceed as in case $(a)$. For $p\ge p_{i_0}$ note that $S^{\phi}_p\in L_{i_0}$, so these points lie on the straight line passing through $S_{p_{i_0}}$ and having slope $\mathbf{a}_{\iota}$. Recall that $L_{i_0}$ coincides with $D_{\mathbf{a}_{\iota}}^{-}:=\lim_{\varepsilon\rightarrow 0}D_{\mathbf{a}_{\iota}-\varepsilon}$ on $[p_{i_0},+\infty)$. Thus \eqref{Regulariz3-equ1} holds also for all such $p$ when $t\in(-\infty,\mathbf{a}_{\iota})$ but the supremum is not attained for any $t$ belonging to this set.
\end{itemize}

\item[$(c)$] Assume that $\phi(T)=+\infty$ for some $T\in\RR$. We distinguish:
\begin{itemize}
\item[$(c.1)$] If there exist infinitely many principal indices, then we can proceed as in case $(a)$ above.

\item[$(c.2)$] If there exist only finitely many principal indices, then the result follows analogously as in case $(b.2)$ before with $\mathbf{a}_{\iota}$ being replaced by $T$.
\end{itemize}
Recall that in case $(c)$ we have \eqref{Regulariz3-equ1} even if $a_p=+\infty$ for all $p\ge 1$ large.
\end{itemize}
\qed\enddemo

\bibliographystyle{plain}
\bibliography{Bibliography}

\end{document}